\documentclass[10pt,3p,times]{elsarticle}  

\usepackage{lineno}

\usepackage{amsmath,amsfonts,amssymb}
\usepackage{amsthm}

\usepackage{etoolbox} 

\newcommand*\linenomathpatch[1]{%
  \cspreto{#1}{\linenomath}%
  \cspreto{#1*}{\linenomath}%
  \csappto{end#1}{\endlinenomath}%
  \csappto{end#1*}{\endlinenomath}%
}

\linenomathpatch{equation}
\linenomathpatch{gather}
\linenomathpatch{multline}
\linenomathpatch{align}
\linenomathpatch{alignat}
\linenomathpatch{flalign}

\usepackage{bm}
\usepackage{graphicx,epstopdf}
\usepackage[position=top,labelfont=normalfont,textfont=normalfont,singlelinecheck=off,justification=raggedright]{subfig}
\usepackage{floatrow}
\usepackage[dvipsnames]{xcolor}
\usepackage{setspace}
\usepackage{enumerate,etaremune}
\usepackage{booktabs}
\usepackage{tabularx,multirow}
\usepackage{framed}
\usepackage{hhline}
\usepackage{url}
\usepackage[unicode,bookmarks=false]{hyperref}
\hypersetup{
  colorlinks,
  citecolor=Blue,linkcolor=Blue,urlcolor=Blue
}
\usepackage[toc,page]{appendix}

\usepackage{tikz}
\usetikzlibrary{shapes.geometric}
\usetikzlibrary{shapes,arrows}
\usetikzlibrary{plotmarks}

\usetikzlibrary{calc}

\usepackage{algorithm}
\usepackage{algorithmicx,algpseudocode}
\usepackage{cases}
\algdef{SE}[DOWHILE]{Do}{doWhile}{\algorithmicdo}[1]{\algorithmicwhile\ #1}

\newcommand{\ie}{{\it i.e.}}

\newcommand{\etal}{{\it et al.}}
\newcommand{\tensor}[1]{\bm{#1}}
\newcommand{\stress}{\sigma}

\newcommand{\tstress}{\tensor{\stress}}

\newcommand{\rn}[1]{\uppercase\expandafter{\romannumeral #1\relax}}

\newsavebox{\dotbox}

\theoremstyle{remark}
\newtheorem{remark}{Remark}

\newcolumntype{L}[1]{>{\raggedright\let\newline\\arraybackslash\hspace{0pt}}m{#1}}
\newcolumntype{C}[1]{>{\centering\let\newline\\arraybackslash\hspace{0pt}}m{#1}}
\newcolumntype{R}[1]{>{\raggedleft\let\newline\\arraybackslash\hspace{0pt}}m{#1}}

\allowdisplaybreaks

\biboptions{sort&compress,square,comma,numbers}

\AtBeginDocument{\hypersetup{citecolor=MidnightBlue,linkcolor=MidnightBlue,urlcolor=MidnightBlue}}

\usepackage{etoolbox}
\makeatletter
\patchcmd{\ps@pprintTitle}
{Preprint submitted to}
{~}
{}{}
\makeatother

\usepackage{xcolor}
\definecolor{darkred}{rgb}{0.55, 0.0, 0.0}
\definecolor{skyblue}{rgb}{0.53, 0.81, 0.92}

\newcommand{\vv}{\mathbf{v}}
\newcommand{\uu}{\mathbf{u}}
\newcommand{\VV}{\mathbf{V}}

\newcommand{\DD}{\mathbf{D}}
\newcommand{\QQ}{\mathbf{Q}}

\newcommand{\xx}{\mathbf{x}}

\newcommand{\XX}{\mathbf{X}}

\newcommand{\FF}{\mathbf{F}}

\newcommand{\TT}{\mathbf{T}}

\newcommand{\NN}{\mathbf{N}}

\newcommand{\nn}{\mathbf{n}}

\renewcommand{\tt}{\mathbf{t}}
\renewcommand{\aa}{\mathbf{a}}
\newcommand{\PP}{\mathbf{P}}
\newcommand{\qq}{\mathbf{q}}

\newcommand{\BB}{\mathbf{B}}

\newcommand{\II}{\mathbf{I}}

\renewcommand{\AA}{\mathbf{A}}

\def\argmin{\mathop{\rm argmin}}

\newcommand{\add}[1]{#1}
\newcommand{\del}[1]{\ignorespaces}

\usepackage{hyperref}
\usepackage{booktabs}

\begin{document}

\begin{frontmatter}

\title{BFEMP: Interpenetration-Free MPM-FEM Coupling with Barrier Contact}

\author[UCLA]{Xuan Li\fnref{equal}\corref{corr}}
\ead{xuanli1@math.ucla.edu}
\author[UCLA,UPENN]{Yu Fang\fnref{equal}}
\author[UCLA]{Minchen Li}
\ead{minchen@math.ucla.edu}
\author[UCLA,UPENN]{Chenfanfu Jiang}
\ead{cffjiang@math.ucla.edu}
\cortext[corr]{Corresponding Author}
\fntext[equal]{These authors contributed equally to this work}

\address[UCLA]{Department of Mathematics, University of California, Los Angeles, United States}
\address[UPENN]{Department of Computer and Information Science, University of Pennsylvania, United States}

\journal{~}

\begin{abstract}
This paper introduces BFEMP, a new approach for \add{monolithically} \del{strongly} coupling the Material Point Method (MPM) with the Finite Element Method (FEM) through \emph{barrier} energy-based particle-mesh frictional contact using a variational time-stepping formulation. The fully implicit time integration of the coupled system is recast into a barrier-augmented unconstrained nonlinear optimization problem. A modified line-search Newton's method is adopted to strictly prevent material points from penetrating the FEM domain, ensuring convergence and feasibility regardless of the time step size or the mesh resolutions. The proposed coupling scheme also reduces to a new approach for imposing separable frictional kinematic boundaries for MPM when all nodal displacements in the FEM domain are prescribed with Dirichlet boundary conditions. Compared to standard implicit time integration, the extra algorithmic components associated with the contact treatment only depend on simple point-segment (or point-triangle in 3D) geometric queries which robustly handle arbitrary FEM mesh boundaries represented with codimension-1 simplices. Experiments and analyses are performed to demonstrate the robustness and accuracy of the proposed method.

\end{abstract}


\end{frontmatter}

\section{Introduction}

The Material Point Method (MPM) \cite{sulsky1995application,de2019material} extends the Particle-In-Cell (PIC) \cite{harlow1962particle} and the Fluid Implicit Particle (FLIP) \cite{brackbill1988flip} methods from fluid dynamics to computational solids. In contrast to the commonly used Total Lagrangian Finite Element Method for elastodynamics \cite{hughes2012finite}, MPM utilizes Lagrangian particles to represent continuum materials and an Eulerian background grid to discretizes the governing equations. Except for recent advancements in Total Lagrangian MPM \cite{de2020total,de2021modelling}, MPM is usually considered to be following an Updated Lagrangian kinematic assumption with particles tracking historical deformation, strain, stress, and other constitutive variables through evolving them with velocity fields. The hybrid Lagrangian-Eulerian perspective combined with the Updated Lagrangian kinematics puts MPM in a very advantageous position in modeling and simulating high-speed, large-deformation, and topologically changing events \cite{zhang2016material}. Having gained a lot of attention in the last two decades, MPM and its variants \cite{bardenhagen2004generalized,zhang2011material,sadeghirad2011convected,gan2018enhancement,hu2018moving}  have been successfully applied in challenging problems including multiphase flows, fracture, contact, adaptivity, free-surface flows, soil-fluid mixture, explosives, and granular media \cite{zhang2008material,homel2017field,nairn2003material,homel2018fracture,tan2002hierarchical,zhang2017incompressible,abe2014material,guilkey2007eulerian,gaume2018dynamic}.

Although MPM has been demonstrated effective on a wide range of materials, many application scenarios favor other discretization schemes due to considerations in efficiency, accuracy, and suitability. Correspondingly, a large number of engineering applications necessitate hybrid or coupled solvers combining MPM with other discretization choices. For example, MPM has been coupled or hybridized with the Discrete Element Method (DEM) for solid-fluid interaction \cite{yang2017combined} and granular media \cite{liu2017multi,jiang2020hybrid,chen2021hybrid}, the Finite Difference Method (FDM) for multiphase saturated soils \cite{higo2010coupled}, and the Smoothed Particle Hydrodynamics (SPH) for solid mechanics \cite{raymond2018strategy}.

Even though MPM itself can be derived as a Galerkin Finite Element Method, it is still more common in computational solid mechanics to use the term ``FEM'' to refer to the standard Total Lagrangian mesh-based FEM discretization. In this sense, MPM is much less developed than FEM and still suffers from unique challenges in aspects such as stability, accuracy, boundary condition enforcement, and numerical fracture \cite{cummins2002implicit,guilkey2003implicit,homel2016controlling}. Therefore, FEM is often more suitable for analyzing hyperelastic structures under small or moderate deformations. Resultingly, MPM-FEM coupling becomes highly desirable in many multi-material simulation tasks or those involving strongly heterogeneous deformations\add{, for example, the blast event simulations involving vehicles \cite{goetz2006blast}}. The coupling between MPM and FEM has been extensively studied over the last decade. A natural way to hybridize the two schemes is to treat FEM vertices as MPM particles and embed them into the MPM grid\ \cite{lian2011femp,jiang2015affine}. \add{The FEM shell formulation can also be embedded into the MPM grid \cite{banerjee2005shell}}. In a similar fashion, EMPFE\ \cite{zhang2006explicit} discretizes the entire domain according to the severity of deformation -- small deformation regions with FEM, while large deformation regions with MPM, and then it embeds the displacement of interface FEM vertices to the MPM grid. Despite its simplicity, additional treatment for eliminating the hourglass mode is necessary due to simple trilinear MPM kernels for the interface computation. Later on, AFEMP\ \cite{lian2012adaptive} extends this idea to support the dynamical conversion from severely distorted FEM elements to MPM particles.  These methods handle material interactions through the grid-based MPM contact and inherit common MPM-based contact characteristics such as the strict nonslip condition between contacting interfaces.

To circumvent these issues, CFEMP\ \cite{lian2011coupling} was proposed to only use the interface MPM grid for contact detection while computing frictional contact forces directly based on the contact conditions. CFEMP has been successfully applied to coupling FEM membranes with MPM solids\ \cite{lian2014coupling}, and also to modeling needle-tissue interactions\ \cite{li2021novel}. However, these grid-based MPM-FEM coupling strategies require the FEM boundary element size to be similar to the MPM grid spacing. If it is too large, interpenetration can happen, while if it is too small, intrinsic damping will appear\ \cite{lian2011femp}, and the time step size for explicit time integration would also be more restricted. Accordingly, Cheon and Kim\ \cite{cheon2018efficient} proposed to add extra distributed interaction (DI) nodes on the FEM boundary elements to improve contact detection for large-sized elements. 

An improvement to CFEMP that also applies to AFEMP was later proposed to couple FEM with MPM by handling contact primitive pairs between FEM boundary elements and nearby MPM particles\ \cite{chen2015improved,wu2018coupled}. A penalty method is applied for computing the frictional contact forces. Later, Song \etal\ \cite{song2020non} extended this idea with an iterative contact force computation approach for the simultaneous satisfaction of all contact conditions, together with an improved local search method to prevent interpenetration issues at the contact crack. Bewick\ \cite{bewick2004combined} proposed to insert intermediate nodes at the interface for 1D impact-resistant design problems, and the coupling forces are calculated by FEM displacements which are determined by MPM particles.

So far, all the discussed MPM-FEM coupling works are designed for explicit time integration. Imposing frictional contact between FEM and MPM within implicit time integration is challenging because the associated inequality constraints that need to be simultaneously enforced while solving the nonlinear system of equations are also nonlinear and non-smooth. 
Aulisa \etal\ \cite{aulisa2019monolithic} proposed a monolithic coupling method for implicit MPM and FEM through a conforming interface mesh, while extra care is needed to avoid the sticky artifact in receding contact cases.
Larese \etal\ \cite{larese2019implicit} discussed implicit MPM-FEM coupling researches in geomechanics. In a soil-structure interaction problem, the impact forces on the interface are transferred between the soil and the structure solver and iterated until convergence. A similar method for enforcing nonconforming boundary conditions for MPM\ \cite{chandra2021nonconforming} was also applied. However, as pointed out in their work, the method requires smaller time increments for problems with high relative velocity towards the boundary, as otherwise, the boundary enforcement will be too late, and the incoming
material points may penetrate the nonconforming boundary
surfaces. 

This paper explores MPM-FEM coupling under the assumption of implicit integration in both domains. Compared to explicit time integration, implicit schemes permit substantially larger time steps with superior stability for stiff nonlinear problems. Implicit MPM/GIMP has been explored by many researchers \cite{cummins2002implicit,guilkey2001implicit,guilkey2003implicit,nair2012implicit,charlton2017igimp,love2006unconditionally}. We refer to these literature for more discussions about the advantages of implicit time integration. Our work follows the variational formulation of a wide family of implicit time integrators \cite{ortiz1999variational,kane2000variational}, \del{using which} \add{where} the displacement evolvement in each time step is formulated as the stationarity point of a time discretized energy functional. The resulting optimization-based time integrator has been applied in MPM  with Newton-Krylov \cite{gast2015optimization}, and more recently, quasi-Newton L-BFGS \cite{wang2019hierarchical} solvers.

In this work, integrating both MPM and FEM domains implicitly, we study MPM-FEM coupling based on contact mechanics (thus we do not consider objects with partially mixed discretization choices). 
A barrier-augmented variational frictional contact formulation is known as the Incremental Potential Contact (IPC) \cite{li2020incremental,li2020robust,eipc} was recently proposed for nonlinear elastodynamics with linear kernel FEM, which has also shown to be effective for codimensional models\ \cite{Li2021CIPC}, reduced space dynamics\ \cite{Ferguson:2021:RigidIPC,Lan2021MIPC}, and embedded interfaces\ \cite{choo2021barrier}. It formulates the contact problem during time stepping as minimizing a potential energy inside the manifold of interpenetration-free displacement trajectories characterized by boundary geometric primitives of elastic structures. Extending this approach, we model MPM-FEM coupling as jointly finding optimal FEM mesh nodal displacements and MPM grid nodal displacements under the constraint that MPM particles, with their trajectories \emph{embedded} in grid nodal displacements, maintain strict positive distances to the FEM mesh throughout the implicit integration. The resulting method is named barrier FEMP (BFEMP) because these constraints are enforced using \emph{barrier} energies.
\del{Even though contact conditions are defined between MPM particles and the FEM mesh, the particles remain in fact embedded quadrature points in the MPM grid degrees of freedom, which are real displacement unknown variables that the implicit time integration solves for.}
\add{Even though contact conditions are defined between MPM particles and the FEM mesh, the real displacement unknown variables for the MPM domain which the implicit time integration solves for are still defined on MPM grids. The MPM particles remain embedded quadrature points in the MPM grid degrees of freedom.}
Compared to soft penalty-based methods such as the particle-to-surface contact algorithm recently proposed by Nakamura \etal \cite{nakamura2021particle}, BFEMP requires no stiffness parameter tuning and guarantees strict, hard non-penetration conditions under convergence. Another useful feature of the proposed coupling scheme is that it enables a new way of imposing irregular, separable, and frictional kinematic boundaries for MPM. Irregular boundaries for MPM is a recently advanced topic \cite{tjung2020modeling}. BFEMP inherently enables it by assigning all nodes in the FEM domain with prescribed Dirichlet displacements and letting MPM interacts with them.
\section{Governing Equations} \label{sec:governing}

In this study we focus on elastodynamics based on continuum mechanics. 
The corresponding governing equations for a deformed continuum domain $\Omega^t$ with $\xx\in \Omega^t$ and time  $t\in[0,\infty)$ are given by~\cite{bonet2008nonlinear}
\begin{align}
\frac{D \rho}{Dt} + \rho \nabla^{\mathbf{x}}\cdot \mathbf{v} &= 0,\\
\rho \frac{D \mathbf{v}}{dt} = \nabla^{\mathbf{x}} \cdot \tstress &+ \rho \mathbf{g}, \label{eqn:eulerian-momentum-balance}
\end{align}
where $\rho(\mathbf{x},t)$ is the density, $\mathbf{v}(\xx,t)$ is the velocity, $\tstress(\xx,t)$ is the Cauchy stress, and $\mathbf{g}$ is the gravitational acceleration which is assumed to be the only body force. Under the finite strain assumption (as we shall assume throughout this paper), deformation $\phi(\mathbf{X},t)$ maps $\mathbf{X}\in \Omega^0$ from the material space to $\mathbf{x} \in \Omega^t$ in the world space: $\xx=\phi(\XX,t)$. The deformation gradient is defined as 
\begin{align}
\mathbf{F} = \frac{\partial \phi}{\partial \mathbf{X}} (\mathbf{X},t)
\label{eq:deform_grad}
\end{align}
to describe the local deformation.

In Lagrangian Finite Element analysis for nonlinear dynamics, it is often preferred to derive the weak form for $\XX\in\Omega^0$. Here we can pull back the momentum equation to the material space. Denoting the first Piola-Kirchhoff stress $\mathbf{P} = \mathbf{P}(\mathbf{X},t)$, the Lagrangian momentum equation is then 
\begin{align}
R_0 \mathbf{A}(\mathbf{X},t) = \nabla^{\mathbf{X}}\cdot\mathbf{P}(\mathbf{X},t) + R_0 \mathbf{g}, \label{eqn:lagrangian-momentum-balance}
\end{align}
where $R_0 = R(\mathbf{X},0)$ is the material density at time 0, and $\mathbf{A}(\mathbf{X},t)=\frac{\partial^2 \Phi}{\partial t^2}(\mathbf{X},t)$ is the Lagrangian acceleration, $\mathbf{V}(\XX,t)=\frac{\partial \Phi}{\partial t}(\mathbf{X},t)$ is the Lagrangian velocity. The Cauchy stress is related to the first Piola Kirchhoff stress as 
\begin{equation}
\tstress = \text{det}(\FF)^{-1} \PP \FF^T.
\end{equation}
In this paper we are concerned with hyperelasticity. Thus there exists an elastic energy density function $\psi(\FF)$ such that
\begin{align}
\PP(\FF) = \frac{\partial \psi}{\partial\FF}(\FF).
\end{align}
Without loss of generality, we focus on isotropic materials and adopt a compressible Neo-Hookean constitutive model with 
\begin{align}
\psi(\FF) &= \frac{\mu}{2}\text{tr}(\FF^T\FF-\II)-\mu\log(J) + \frac{\lambda}{2}\log(J)^2,
\end{align}
where $J=\det{\FF}$, $\mu$ and $\lambda$ are the Lam\'e parameters.

\begin{remark}
Many hyperelasticity models such as the Neo-Hookean model are only well defined for $J>0$. Discretely this will impose a nonlinear strict inequality constraint on the displacements for each quadrature point with a discrete sample of $\FF$. Ignoring these constraints in a nonlinear optimization-based Newton solver will cause floating-point number failures when a search step tries to evaluate energy or stress quantities at an intermediate configuration with $J\le 0$. Instead of requiring a reduction of the time step, we directly enforce this constraint through a line search filtering strategy (see Section \ref{sec:optimization}).
\end{remark}

\subsection{Weak Form}

Given trial function $\QQ(\cdot,t):\Omega^0\rightarrow\mathbb{R}^{3}$, the corresponding \emph{Lagrangian} weak form of Equation (\ref{eqn:lagrangian-momentum-balance}) is 
\begin{align}
\int_{\Omega^0} Q_{\alpha}(\XX,t) R_0 A_{\alpha}(\XX,t)d\XX = \int_{\partial\Omega^0} Q_{\alpha}T_{\alpha} dS(\XX) - \int_{\Omega^0} Q_{\alpha,\beta} P_{\alpha\beta}d\XX + \int_{\Omega^0} Q_{\alpha}(\XX,t) R_0 g_{\alpha} d\XX, \label{eqn:lagrangian-weak-form}
\end{align}
where $T_{\alpha} = P_{\alpha\beta}N_{\beta}$ (with $\NN(\XX)$ being the material space normal) is the traction field at the domain boundary $\partial\Omega^0$, on which one could presribe traction boundary conditions as needed.

While Total Langrangian FEM typically discretizes Equation (\ref{eqn:lagrangian-weak-form}), MPM usually adopts the Updated Lagrangian view and consequently discretizes an \emph{Eulerian} weak form instead \cite{zhang2016material}. Correspondingly the stress derivatives are discretized at the current configuration $\Omega^t$. We can either push forward Equation (\ref{eqn:lagrangian-weak-form}) or directly integrate Equation (\ref{eqn:eulerian-momentum-balance}) to reach
\begin{align}
\int_{\Omega^t} q_{\alpha}(\xx,t)\rho(\xx,t)a_i(\xx,t)d\xx &= \int_{\partial\Omega^t}q_{\alpha}t_{\alpha}ds(\xx) - \int_{\Omega^t}q_{\alpha,\beta}\sigma_{\alpha\beta}d\xx + \int_{\Omega^t}q_{\alpha}(\xx,t)\rho(\xx,t)g_{\alpha}d\xx,
\end{align}
where $\qq(\xx,t)=\QQ(\Phi^{-1}(\xx,t),t)$ is the push forward of $\QQ$, \ie, an Eulerian trial function,   $\aa(\xx,t)=\AA(\Phi^{-1}(\xx,t),t)=\frac{D\vv}{Dt}(\xx,t)=\frac{\partial\vv}{\partial t}+\vv\cdot\nabla \vv$, and $\tt = \tstress \nn$ is the traction at $\partial \Omega^t$ with $\nn(\xx)$ being the outward pointing normal.

\subsection{Incremental Variational Form}

To enable the development of efficient optimization-based time integrators, we follow the variational treatment of the time-discretized incremental problem \cite{ortiz1999variational,radovitzky1999error}. Concretely, for a broad family of time discretization schemes (we will focus on backward Euler and the Newmark-$\beta$ family in this paper), the solution of Equation $(\ref{eqn:lagrangian-weak-form})$ during an incremental time step, \ie, the advancement from $\phi^n=\phi(\XX,t^{n})$ to $\phi^{n+1}=\phi(\XX,t^{n+1})$), for an hyperelastic material is given by minimizing the following functional:
\begin{align}
I(\phi^{n+1}) &= \int_{\Omega_0} \left(\frac{1}{2} \frac{R_0}{\beta \Delta t^2}\|\phi^{n+1}\|^2 + 2\alpha \Psi(\FF^{n+1})\right)d\XX - \int_{\Omega_0} R_0 \bar{\BB}_{n+1}\cdot\phi^{n+1}d\XX - 2\alpha\int_{\partial\Omega_0}\TT\cdot\phi^{n+1}dS(\XX), \label{eqn:continuous-IP}
\end{align}
where 
\begin{align}
\bar{\BB}_{n+1} = 2\alpha\mathbf{g} + \frac{1}{\beta \Delta t^2} \left( \phi^n + \Delta t \frac{\partial\phi}{\partial t}(\XX,t^n) +\alpha\left(1-2\beta\right)\Delta t^2  \frac{\partial^2\phi}{\partial t^2}(\XX,t^n)\right)
\end{align}
encodes the inertia term and is a constant field in the minimization problem. The velocity update rule is
\begin{align}
    \frac{\partial\phi}{\partial t}(\XX,t^{n+1}) = \frac{\partial\phi}{\partial t}(\XX,t^{n}) + \Delta t \left(\left(1-\gamma\right)\frac{\partial^2\phi}{\partial t^2}(\XX,t^n) + \gamma \frac{\partial^2\phi}{\partial t^2}(\XX,t^{n+1}) \right).
\end{align}
Here we have slightly modified the energy proposed by Radovitzky and Ortiz \cite{radovitzky1999error} to let it be compatible with backward Euler since their version was explicitly built for Newmark's algorithm.

Note that in the case of backward Euler ($\alpha = 1$, $\beta = \frac{1}{2}$, $\gamma = 1$), it can be easily shown that the Euler-Lagrangian equation of functional (\ref{eqn:continuous-IP}) gives back the time-discretized momentum balance
\begin{align}
\frac{R_0}{\Delta t^2} \phi^{n+1} -\nabla^{\XX}\cdot\PP^{n+1} = R_0 \mathbf{g} + \frac{R_0}{\Delta t^2}(\phi^n + \Delta t \VV^n).
\end{align}
Similarly, for Middle-point Newmark ($\alpha = \frac{1}{2}$, $\beta = \frac{1}{4}$, $\gamma = \frac{1}{2}$), we would get 
\begin{align}
\frac{R_0}{\Delta t^2} \phi^{n+1} -\frac{1}{4}\nabla^{\XX}\cdot\PP^{n+1} = \frac{1}{4}R_0 \mathbf{g} + \frac{R_0}{\Delta t^2}(\phi^n + \Delta t \VV^n + \frac{1}{4}\Delta t^2 \AA^n).
\end{align}
Both of them match the results from temporally discretizing the Lagrangian form of Equation (\ref{eqn:eulerian-momentum-balance}).

\begin{figure}[ht]
    \centering
    \includegraphics[width=0.9\linewidth]{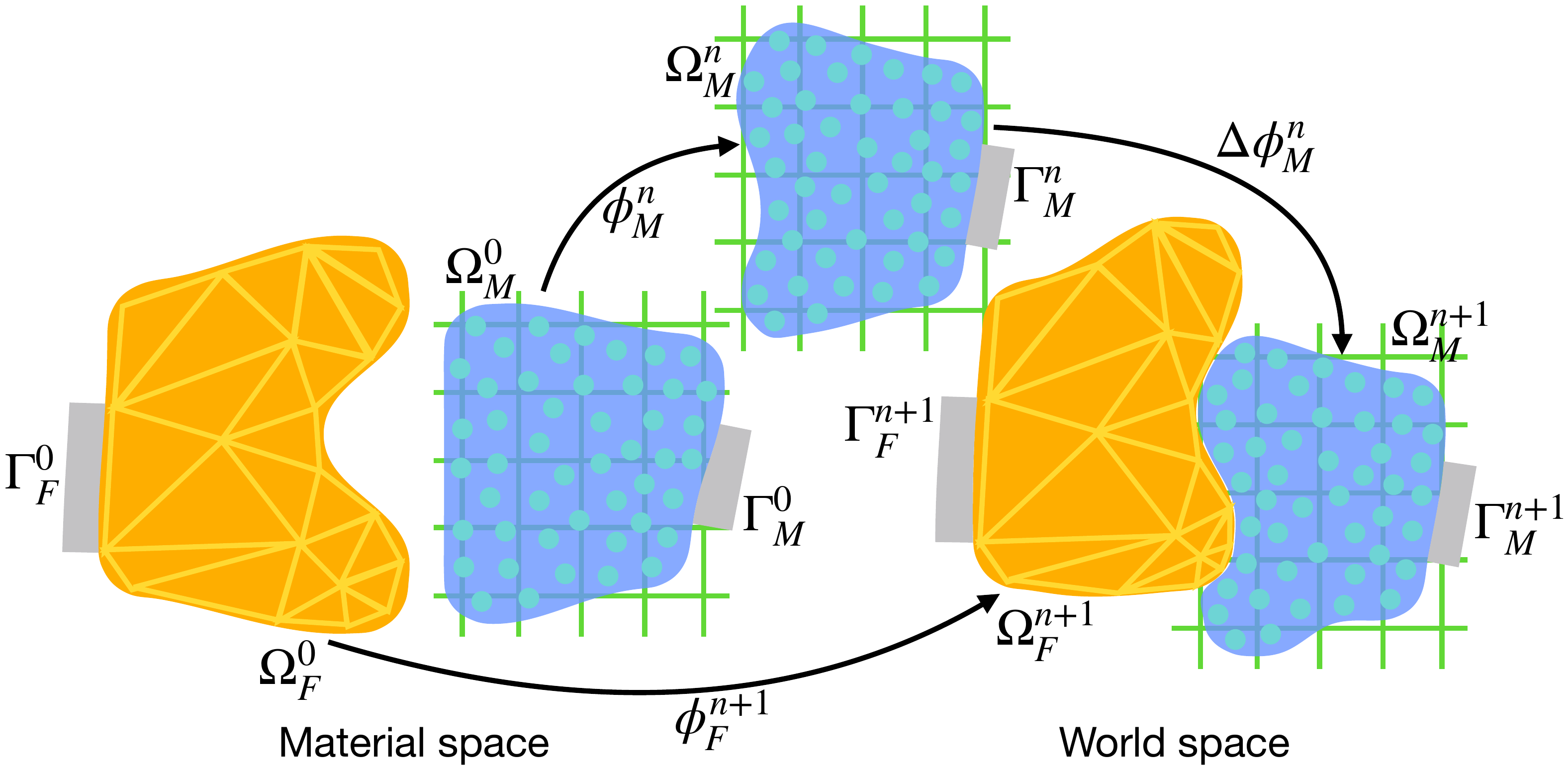}
    \caption{\textbf{Deformation map.} The material space of FEM and MPM domains ($\Omega_\text{F}^0$ and $\Omega_\text{M}^0$) on the left are mapped via $\phi_\text{F}$ and $\phi_\text{M}$ to the world space ($\Omega_\text{F}^{n+1}$ and $\Omega_\text{M}^{n+1}$) on the right, with $\Omega_\text{F}^t \cap \Omega_\text{M}^t = \varnothing$ for all $t\in [0,\infty)$. With Updated Lagrangian kinematics, MPM treats $\Omega_M^n$ as the ``\add{intermediate}\del{intemediate}'' material space and focuses on studying the deformation from $\Omega_M^n$ to $\Omega_M^{n+1}$. Here $\Gamma$ represents Dirichlet boundary or nonzero Neumann boundary.}
    \label{fig:map}
\end{figure}

\subsection{Coupling}
In the proposed framework, the Lagrangian form is discretized on the FEM domain $\Omega_F$ (Section\ \ref{sec:discretize_FEM}), while the Eulerian form is discretized on the MPM domain $\Omega_M$ (Section\ \ref{sec:discretize_MPM}). We assume 
\begin{align}
\Omega_F^0 \cap \Omega_M^0 = \varnothing \label{eqn:init-no-intersection}
\end{align}
and their corresponding deformation map $\xx_M = \phi_M(\XX_M,t)$ and $\xx_F=\phi_F(\XX_F,t)$ are fully governed by their own variational forms in the absence of any coupling mechanism. In other words, without coupling, evolving the two domains independently (with two solvers) is equivalent to minimizing 
\begin{align}
\Pi(\phi_M,\phi_F) = I(\phi_M) + I(\phi_F) \label{eq:sum-ip}
\end{align}
by directly letting $\Omega=\Omega_F \cup \Omega_M$.

The coupling between FEM and MPM domains are then modeled via imposing the non-interpenetration constraints:
\begin{align}
\Omega_F^t \cap \Omega_M^t = \varnothing, \quad \forall t\in[0,\infty) \label{eq:set-intersection}
\end{align}
between the two domains. Note that \add{the time is discretized with} \del{when discretizes in time,} $t=t^0, t^1, t^2, \dots, t^n$. We use $\Omega^n$ to denote $\Omega^{t^n}$ for notational simplicity.  See Figure \ref{fig:map} for an illustration. Thus the final variational MPM-FEM coupling problem can be described as minimizing Equation (\ref{eq:sum-ip}) under the equality constraint described by Equation (\ref{eq:set-intersection}).

The feasible region described by Equation (\ref{eq:set-intersection}) can be equivalently expressed as
\begin{align}
\phi_M(\XX_M, t) \neq \phi_F(\XX_F,t), \quad \forall \XX_M\in\Omega_M^0, \XX_F\in\Omega_F^0, t\in[0,\infty). \label{eqn:phi-neq}
\end{align}
Moreover, if we define
\begin{align}
d(\phi_M,\phi_F,t) = \min_{\XX_F,\XX_M} \|\phi_F(\XX_F,t)-\phi_M(\XX_M,t)\|
\end{align}
to describe the Euclidean proximity between $\Omega_M^t$ and $\Omega_F^t$, Equation (\ref{eqn:phi-neq}) can be further converted to a strict inequality constraint
\begin{align}
d(\phi_M,\phi_F,t) > 0, \quad \forall t\in[0,\infty).
\label{eq:min_d_contact_form}
\end{align}
Note that in Equation (\ref{eqn:init-no-intersection}), we have assumed that the undeformed domains are non-overlapping. Thus the minimization problem starts with a strictly feasible solution at $t=0$.

In Section \ref{sec:contact-potential} we describe a barrier method that results in a contact pressure for enforcing Equation (\ref{eqn:init-no-intersection}). See Section \ref{sec:friction} for extra components on \add{incorporating} \del{incooporating} tangential frictional effects.
\section{Discretization}\label{sec:discretization}
The finite element domain is discretized with linear simplex elements (triangles in 2D and tetrahedra in 3D), and the material point domain is discretized using a collection of material points and an Eulerian background grid with quadratic B-spline kernels. Both schemes adopt mass lumping and assume zero traction at boundaries unless otherwise specified in an example. 
Stacking all nodal positions, velocities, and accelerations from both the FEM mesh and the MPM grid at time step $n$ as $x^n = [(x^n_\text{F})^T, (x^n_\text{M})^T]^T$, $v^n = [(v^n_\text{F})^T, (v^n_\text{M})^T]^T$, and $a^n = [(a^n_\text{F})^T, (a^n_\text{M})^T]^T$  where subscripts F stands for FEM while M for MPM, the unified time integration update rule can be written as
\begin{equation}
    \begin{split}
        \tilde{v}^{n+1} &= v^n + \Delta t((1-\gamma) a^{n} + \gamma a^{n+1}), \\
        \tilde{x}^{n+1} &= x^n + \Delta t v^n + \alpha \Delta t^2 ((1-2\beta) a^{n} + 2\beta a^{n+1})),
    \end{split}
    \label{eq:discrete_update_rule}
\end{equation}
and it is equivalent to first solving the optimization problem
\begin{equation}
    \min_{x}:\left( \frac{1}{2} \|x - \hat{x}^{n}\|^2_M + 2 \alpha \beta \Delta t^2 \Psi(x)\right)
    \label{eq:discrete_IP}
\end{equation}
to get $\tilde{x}^{n+1}$, and then explicitly calculating $\tilde{v}^{n+1}$. Here $\hat{x}^{n} = x^n + v^n \Delta t +  \alpha  (1-2\beta) \Delta t^2 a^n $, and $\Psi(x) = \sum_q V_q^0 \psi(\FF_q(x))$ with $q$ belonging to FEM elements or MPM particles and $V_q^0$ the rest volume of a FEM element or a MPM particle. 

For FEM, the time integration is solely performed on the Lagrangian nodal degrees of freedom throughout the simulation and so $\tilde{x}^{n+1}_\text{F} = x^{n+1}_\text{F}$ and $\tilde{v}^{n+1}_\text{F} = v^{n+1}_\text{F}$. But for MPM, Equations (\ref{eq:discrete_update_rule}) are only part of the Eulerian time integration performed on the Eulerian grid before and after particle-grid transfers for $x_\text{M}$ and $v_\text{M}$, so $\tilde{x}^{n+1}_\text{M} \neq x^{n+1}_\text{M}$ and $\tilde{v}^{n+1}_\text{M} \neq v^{n+1}_\text{M}$ (see Section\ \ref{sec:discretize_MPM} for details).
Note that the minimization problem (\ref{eq:discrete_IP}) is equivalent to the discrete form of the variational time integration in\ \cite{ortiz1999variational,li2019decomposed,wang2019hierarchical} for hyperelastic problems. $x_\text{F}$ and $x_\text{M}$ are coupled through contact modeling between the two domains (Section\ \ref{sec:contact}).

\subsection{The Finite Element Domain}
\label{sec:discretize_FEM}

For the FEM domain, nodal positions and velocities are stored on mesh vertices and updated directly. The nodal \add{masses} \del{mass} are kept constant, \ie, $m^n_i = m_i$.

Inside any simplex element, the material and world space coordinate of an arbitrary location $\XX$ are expressed using linear interpolation kernels $N_i(\XX)$ of node $\XX_i$ as
\begin{equation}
    \XX = \sum_i N_i(\XX) \XX_i \quad \text{and} \quad \phi(\XX,t) = \sum_i N_i(\XX) \phi(\XX_i,t).
\end{equation}
Then according to the definition of $\FF$ (Equation (\ref{eq:deform_grad})), with $N_i$ being the linear hat function, the deformation gradient is piecewise constant. Inside a linear simplex element $e$ it is directly evaluated as a function of $x$:
\begin{equation}
    \FF_e(x) = \TT_e(x) \BB_e^{-1},
\end{equation}
where $\TT_e(x)$ is the current triangle basis of element $e$ and $\BB_e$ is the triangle basis of element $e$ in material space.

\subsection{The Material Point Domain}
\label{sec:discretize_MPM}

For the MPM domain, the nodal positions $x^n_\text{M}$ are the uniform Cartesian grid coordinates at each time step. The grid velocity $\vv_i^n$ and mass $m_i^n$ are transferred from particles. The nodal movements are conceptual. $\tilde{x}^{n+1}_\text{M}$ and $\tilde{v}^{n+1}_\text{M}$ will be transferred back to particles for advection.

Similar to FEM, each MPM grid node $i$ is associated with a kernel function $N_i(\xx)$ for the grid to represent the continuous field. Note that the kernel is defined in terms of $\xx$ rather than $\XX$ because the grid is essentially a discretization of $\Omega_M^n$ -- a direct consequence of adopting Updated Lagrangian kinematics. When $N_i$ is evaluated at a particle location $\xx_q^n$, a shorter notation $N_i(\xx_q^n) = w_{iq}^n$ is from now on used instead. 
Here $N_i$ directly takes the current particle location $\xx_q^n$ as input as opposed to FEM because there is no globally defined reference configuration in MPM and the deformation is evolved over time steps rather than recomputed using a rest shape. More specifically, the deformation gradient of a particle $q$ is defined as
\begin{equation}
    \FF_q(x) = \sum_i \xx_i (\nabla w^n_{ip})^T \FF_q^n.
\end{equation}
In this paper, without loss of generality, we adopt the quadratic B-spline kernel for $N_i(\xx)$ to avoid MPM's cell-crossing instability \cite{steffen2008analysis}. Other kernels based on the NURBS \cite{long2019using}, Generalized Interpolation Material Point Method (GIMP) \cite{bardenhagen2004generalized}, Convective Particle Domain Interpolation (CPDI) \cite{homel2016controlling,sadeghirad2013second,nguyen2017family}, or the Dual Domain Material Point (DDMP) \cite{long2016representing,zhang2013dual} can also be directly used in our framework.

To transfer information between the particles and the grid, we implemented options including the Affine Particle-In-Cell (APIC) method \cite{jiang2017angular,jiang2015affine} (Table\ \ref{tb:apic}), Particle-In-Cell (PIC) method\ \cite{harlow1962particle} (Table\ \ref{tb:pic}), and the Fluid-Implicit Particle (FLIP) method \cite{brackbill1988flip} (Table\ \ref{tb:flip}). Note that other transfer schemes such as XPIC \cite{hammerquist2017new} can also be applied in our framework in a straightforward manner.

\begin{table}[h!]
\caption{APIC Particle-Grid Transfer}
\vspace{2pt}
\centering
\begin{tabular}{l|l}
\toprule
Particles to grid (APIC) & Grid to particles (APIC) \\ \hline
\parbox{7cm}{
\begin{equation*}
\begin{split}
       m^n_i &= \sum_p m_p w_{ip}^n\\
       \DD_p^n &= \sum_i w_{ip}^n (\xx_i^n - \xx_p^n)(\xx_i^n - \xx_p^n)^T\\
       m_i^n \vv_i^n &= \sum_p w_{ip} m_p^n(\vv_p^n + \BB_p(\DD_p)^{-1}(\xx_i^n - \xx_p^n))
\end{split}
\end{equation*}
} &  
\parbox{7cm}{
\begin{equation*}
\begin{split}
       \vv_p^{n+1} =& \sum_i \tilde{\vv}^{n+1}_i w_{ip}^n\\
       \xx_p^{n+1} =& \sum_i \tilde{\xx}^{n+1}_i w_{ip}^n\\
       \BB_p^n =& \frac{1}{2} \sum_i w_{ip}^n \bigg{(}\tilde{\vv}^{n+1}_i(\xx_i^n - \xx_p^n + \tilde{\xx}_i^{n+1} - \xx_p^{n+1})^T \\
       &+ (\xx_i^n - \xx_p^n - \tilde{\xx}_i^{n+1} + \xx_p^{n+1})(\tilde{\vv}^{n+1}_i)^T\bigg{)}\\
       \FF_p^{n+1} =&  \sum_i \tilde{\xx}^{n+1}_i (\nabla w^n_{ip})^T \FF_p^n
\end{split}
\end{equation*}
}
\\ \bottomrule
\end{tabular}
\label{tb:apic}
\end{table}

\begin{table}[h!]
\caption{PIC Particle-Grid Transfer}
\vspace{2pt}
\centering
\begin{tabular}{l|l}
\toprule
Particles to grid (PIC) & Grid to particles (PIC) \\ \hline
\parbox{7cm}{
\begin{equation*}
\begin{split}
       m^n_i &= \sum_p m_p w_{ip}^n\\
       m_i^n \vv_i^n &= \sum_p w_{ip}^n m_p^n \vv_p^n
\end{split}
\end{equation*}
} &  
\parbox{7cm}{
\begin{equation*}
\begin{split}
       \xx_p^{n+1} =& \sum_i \tilde{\xx}^{n+1}_i w_{ip}^n\\
       \vv_{p}^{n+1} =& \sum_{i}  \tilde{\vv}^{n+1}_i w_{ip}^n\\
        \FF_p^{n+1} =&  \sum_i \tilde{\xx}^{n+1}_i (\nabla w^n_{ip})^T\FF_p^n 
        \end{split}
\label{tb:pic}
\end{equation*}
}
\\ \bottomrule
\end{tabular}
\end{table}

\begin{table}[h!]
\caption{FLIP Particle-Grid Transfer}
\vspace{2pt}
\centering
\begin{tabular}{l|l}
\toprule
Particles to grid (FLIP) & Grid to particles (FLIP) \\ \hline
\parbox{7cm}{
\begin{equation*}
\begin{split}
       m^n_i &= \sum_p m_p w_{ip}^n\\
       m_i^n \vv_i^n &= \sum_p w_{ip}^n m_p^n \vv_p^n
\end{split}
\end{equation*}
} &  
\parbox{7cm}{
\begin{equation*}
\begin{split}
       \xx_p^{n+1} =& \sum_i \tilde{\xx}^{n+1}_i w_{ip}^n\\
       \vv_{p}^{n+1} =& \vv_p^n + \sum_{i} w_{ip}^n (\tilde{\vv}^{n+1}_i - \vv^{n}_i)\\
        \FF_p^{n+1} =&  \sum_i \tilde{\xx}^{n+1}_i (\nabla w^n_{ip})^T\FF_p^n 
        \end{split}
\label{tb:flip}
\end{equation*}
}
\\ \bottomrule
\end{tabular}
\end{table}

\section{The Contact between Domains}
\label{sec:contact}

\subsection{Contact Potential} \label{sec:contact-potential}

Recently for Lagrangian FEM, Li \etal \cite{li2020incremental,eipc} proposed a consistent variational contact model that smoothly approximates the nonsmooth contact phenomena with bounded error, and demonstrated its convergence under refinement for piecewise linear boundary discretization. Here we customize the FEM contact potential \cite{eipc} to the FEM-MPM coupling setting by defining the inter-surface contact potential between surfaces $\partial \Omega_\text{M}$ and $\partial \Omega_\text{F}$ to be
\begin{equation}
    \int_{\partial \Omega_\text{M}^0} b(\min_{\xx_f \in \partial \Omega_\text{F}^t} d^{PP}(\xx_m,\xx_f), \hat{d}) \mathbf{d}\XX_m,
    \label{eq:ipc_smooth}
\end{equation}
where $d^{PP}(\xx_m, \xx_f) = \|\xx_m - \xx_f\|$ is the point-point distance function, and
\begin{equation} 
    b(d, \hat{d}) = \begin{cases}
        -\kappa (\frac{d}{\hat{d}} - 1)^2 \ln{(\frac{d}{\hat{d}})} & 0 < d < \hat{d} \\
        0 & d \geq \hat{d}
    \end{cases}
    \label{eq:barrier}
\end{equation}
is a smoothly clamped barrier function that serves as the contact energy density with $d$ the input distance, $\hat{d}$ a small distance threshold below which contact activates, and $\kappa$ in $Pa$ the barrier stiffness\ \cite{li2020incremental,eipc}\add{, which scales the magnitude of contact forces at a certain distance}.

\begin{figure}[h!]
    \centering
    \includegraphics[width=0.4\linewidth]{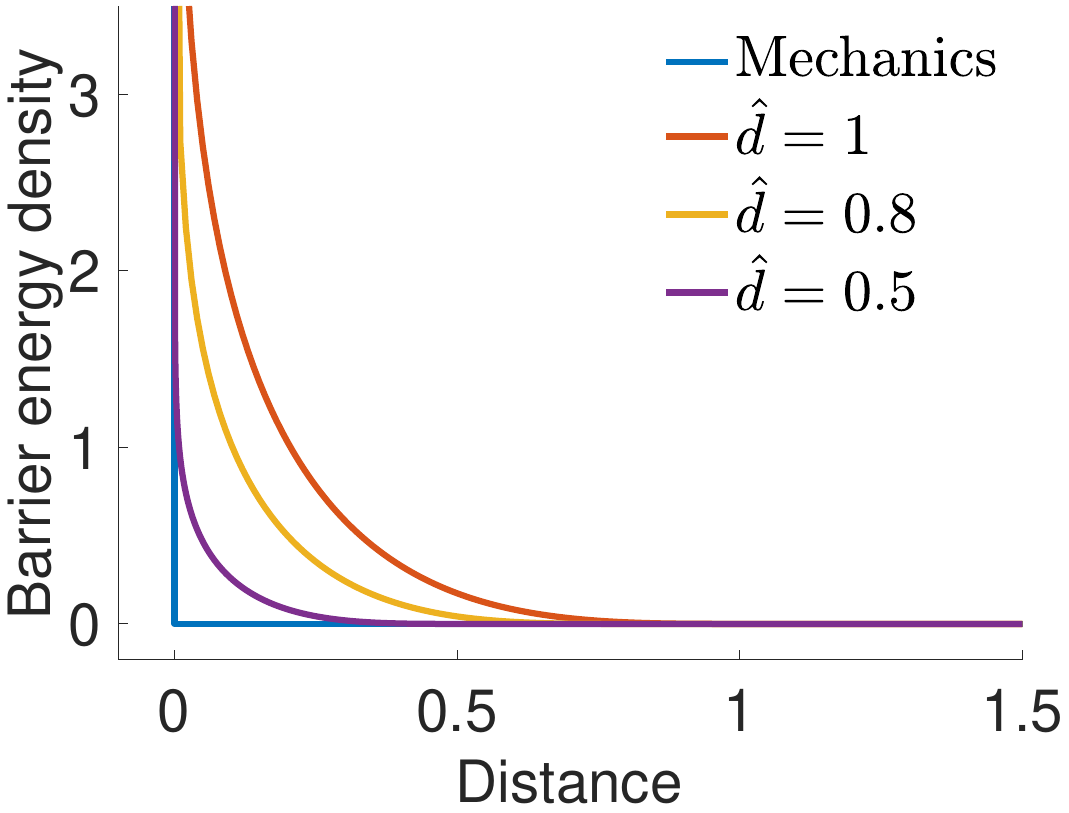}
    \caption{\textbf{The barrier energy density function Equation (\ref{eq:barrier}) plotted with different $\hat{d}$.} Decreasing $\hat{d}$ asymptotically matches the discontinuous definition of the contact condition.}
    \label{fig:barrier_plot}
\end{figure}

Intuitively, $\min_{\xx_f \in \partial \Omega_\text{F}^t} d^{PP}(\xx_m,\xx_f)$ is the distance between a material point $\xx_m\in \partial \Omega_\text{M}^t$ and surface $\partial \Omega_\text{F}^t$, and Equation (\ref{eq:ipc_smooth}) can be viewed as an integration of the point($\xx_m$)-surface($\partial \Omega_\text{F}^t$) contact energy density over surface $\partial \Omega_\text{M}^0$.
The barrier function $b$ smoothly increases from $0$ to infinity as the input distance decreases from $\hat{d}$ to $0$, providing arbitrarily large repulsion to ensure no interpenetration and at the same time bound the contact gap error within $\hat{d}$ (Figure\ \ref{fig:barrier_plot}).
As $\hat{d}\rightarrow 0$, the approximation error between the contact energy density function $b$ and the real contact phenomenon described in Equation (\ref{eq:min_d_contact_form}) decreases, which also makes $\partial \Omega_\text{M}$ and $\partial \Omega_\text{F}$ interchangeable in the limit.
Note that the integration is performed in the material space while the distance is evaluated in the world space.

\subsection{Discretization}
\begin{figure}[t]
    \centering
    \includegraphics[width=0.8\linewidth]{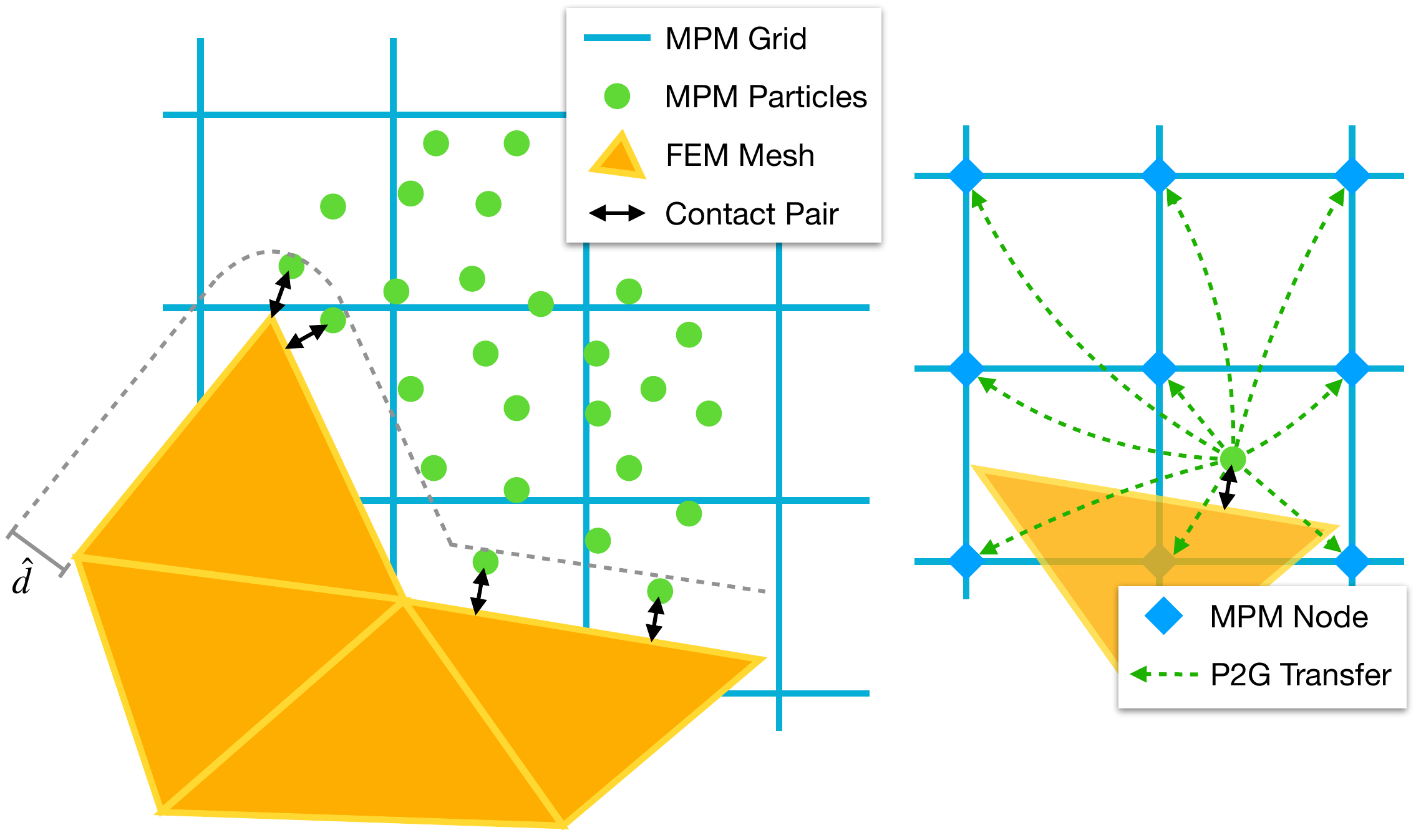}
    \caption{\textbf{Contact constraint pairs.} Left: Contact activates on all pairs of MPM particles and FEM boundary elements with distance below $\hat{d}$. Right: Contact force is transferred from MPM particles to MPM grids via chain rule.}
    \label{fig:constraint}
\end{figure}

After applying FEM and MPM discretization schemes, assuming 2D, the contact potential Equation (\ref{eq:ipc_smooth}) is discretized to be 
\begin{equation}
    B(x) = \sum_{q\in \mathcal{Q}} \omega_q b(\min_{e \in \mathcal{B}} d^{PE}(\xx_q, e), \hat{d}), \label{eq:discrete-potential}
\end{equation}
where $\mathcal{Q}$ is the set of all MPM particles, $\omega_q$ is the integration weight (equivalently, the boundary area) of MPM particle $q$, $\mathcal{B}$ is the set of FEM boundary edges, and $d^{PE}(\xx_q, e)$ is the point-edge distance between particle $\xx_q$ and edge $e$. Note that the MPM grid nodal positions $\tilde{x}_\text{M}$ to be solved and the particle positions $\xx_q$ after advection using $\tilde{x}_\text{M}$ are linearly related through the particle-grid transfer kernel (Figure\ \ref{fig:constraint} right), and so the contact force on the MPM nodal degrees of freedom can be calculated by applying the chain rule:
\begin{equation}
    \frac{\partial B}{\partial x_\text{M}} = \sum_{q\in\mathcal{Q}} \left( \frac{\partial \xx_q}{\partial x_\text{M}} \right)^T \frac{\partial B}{\partial \xx_q},
\end{equation}
where $$\left(\frac{\partial \xx_q}{\partial x_M}\right)_{\alpha,id+\beta} = \delta_{\alpha\beta}w_{iq}^n$$ with $\alpha,\beta=1,2,...,d$ and $d=2$ or $3$ the spatial dimension since $\xx_q = \sum_i w_{iq}^n\xx_i$.

Ideally, $\omega_q$ should be zero for interior particles. It should reveal the proportional surface area of boundary particles. Since FEM boundary elements are always outside the MPM domain and the barrier function $b$ is only activated at a small distance, the activation of $b$ can be applied to conveniently decide whether an MPM particle is at the boundary or the interior without explicitly identifying the MPM domain boundary in each time step (Figure\ \ref{fig:constraint} left). Then assuming a close to uniform particle distribution to be maintained throughout the simulation, $\omega_q$ can be set to $2\sqrt{V_q^0/\pi}$ in 2D and $\pi\left(3V_q^0/\left(4\pi\right)\right)^{\frac{2}{3}}$ in 3D for all particles, which is the area of the largest cross section of a spherical particle $q$.

\begin{figure}[t]
    \centering
    \includegraphics[width=\linewidth]{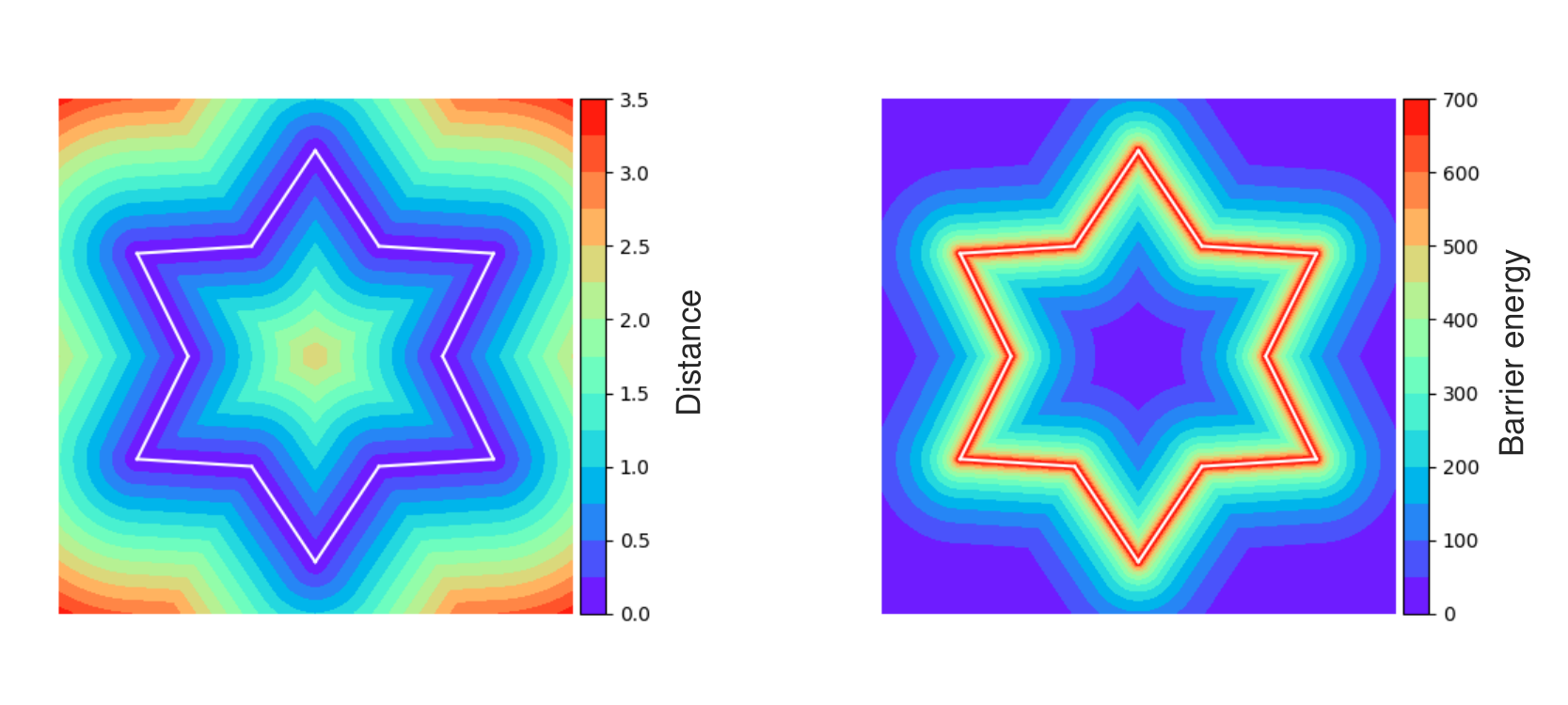}
    \caption{A visual demonstration in 2D of (left) the unsigned distance function to a  segmented mesh, and (right) the corresponding barrier function  (\ref{eq:barrier}) visualized with an exaggerated $\hat{d}$ parameter.}
    \label{fig:distance}
\end{figure}

The minimization operator in the potential (Equation (\ref{eq:discrete-potential})) helps to compute the point-polyline distance from point-edge distances. As the barrier function $b$ is monotonically decreasing, the potential can be rewritten as
\begin{equation}
    B(x) = \sum_{q\in\mathcal{Q}} \omega_q \max_{e \in \mathcal{B}} b( d^{PE}(\xx_q, e), \hat{d}).
\end{equation}
Due to the existence of a $\max$ operator, it is only $C^0$ continuous, making the incremental potential challenging to be efficiently minimized by gradient-based optimization methods like Newton's method. Since the barrier function $b$ is with local support around each boundary element, and that it maps a majority of the activate distances to tiny potential values (Figure\ \ref{fig:distance}), the maximization of the potential field can be well approximated by summation, with the duplicate potential around FEM boundary nodes compensated by subtraction as proposed by Li \etal\ \cite{eipc}:
\begin{equation}
    B(x) = \sum_{q\in\mathcal{Q}} \omega_q \left( \sum_{e \in \mathcal{B}} b( d^{PE}(\xx_q, e), \hat{d}) - \sum_{k \in \hat{\mathcal{B}}} \left(\eta_k-1\right) b(d^{PP}(\xx_q, \xx_k), \hat{d}) \right),
\end{equation}
where $\hat{\mathcal{B}}$ is the set of all FEM boundary nodes and $\eta_k$ is the number of FEM boundary edges incident to node $k$. For closed manifold domains in 2D, $\eta_k = 2$ for all $k$.

Similarly, in 3D, the discretized contact potential becomes
\begin{equation}
B(x) = \sum_{q\in\mathcal{Q}} \omega_q \left( \sum_{t\in \mathcal{B}} b( d^{PT}(\xx_q, t),\hat{d}) - \sum_{e\in \hat{\mathcal{B}}} \left(\eta_e - 1\right) b(d^{PE}(\xx_q, e), \hat{d}) + \sum_{\xx_p\in \tilde{\mathcal{B}}} b(d^{PP}(\xx_q,\xx_p), \hat{d}) \right),
\end{equation}
where $\mathcal{B}$ is now the set of all FEM boundary triangles, $\hat{\mathcal{B}}$ is the set of all edges on the FEM boundary with $\eta_e$ the number of FEM boundary triangles incident to edge $e$, $\tilde{\mathcal{B}}$ the set of all nodes that are in the interior of the FEM boundary surface mesh, and $d^{PT}(\xx_q, t)$ the point-triangle distance between particle $\xx_q$ and triangle $t$. For closed manifold domains in 3D, $\eta_e = 2$ for all $e$.

Adding the contact potential into the incremental potential, the minimization problem for time integration is now fully unconstrained:
\begin{equation}
    \min_{x}: \frac{1}{2} \|x - \hat{x}^{n}\|^2_M + 2 \alpha \beta \Delta t^2 \left( \Psi(x) + B(x) \right).
    \label{eq:discretized_IP_withB}
\end{equation}
Since the distance values measured for contacting MPM particle - FEM simplex pairs are all unsigned, Problem (\ref{eq:discretized_IP_withB}) may contain a local optimum at configurations with intersections. To be consistent with the continuous constraint Equation (\ref{eq:min_d_contact_form}), it is also constrained that the iterates always stay in the feasible region on one side of the barrier without crossing. This is achieved by applying the interior-point filter line-search algorithm\ \cite{wachter2006implementation} with continuous collision detection (CCD)\ \cite{brochu2012efficient,Li2021CIPC}.

\subsection{Friction}\label{sec:friction}
To model frictional contact, local frictional forces $F_k$ can be added for every active contact pair $k$.
For each such pair $k$, at the current state $x$, a consistently oriented sliding basis $T_k(x) \in \mathbb{R}^{d m\times(d-1)}$ can be constructed, where $m$ is the total number of colliding nodes and $d$ is the dimension of space, such that $\uu_k = T_k(x)^T (\Delta t v^n + \Delta t^2 ((1-\gamma) a^{n} + \gamma a^{n+1})) \in \mathbb{R}^{d-1}$ provides the local relative sliding displacement in the frame orthogonal to the distance gradient. The corresponding sliding velocity is then $\vv_k = \uu_k/\Delta t \in \mathbb{R}^{d-1}$.

Maximizing dissipation rate subject to the Coulomb constraint defines friction forces variationally\ \cite{moreau2011unilateral,goyal1991planar}
\begin{align}
\begin{split}
    F_k(x) = T_k(x) \> \argmin_{\boldsymbol{\beta} \in \mathbb{R}^{d-1}} \boldsymbol{\beta}^T \vv_k \quad \text{s.t.} \quad  \|\boldsymbol{\beta}\| \leq \mu \lambda_k,
   \end{split}
 \label{eq:frictionForceDef}
\end{align}
where $\lambda_k$ is the contact force magnitude 
and $\mu$ is the local friction coefficient. This is equivalent to
\begin{equation}
    F_k(x) = - \mu \lambda_k T_k(x) f(\|\uu_k\|) \> \mathbf{s}(\uu_k),
    \label{eq:smooth-frict}
\end{equation}
with $\mathbf{s}(\uu_k) = \frac{\uu_k}{\|\uu_k\|}$ when $\|\uu_k\|> 0$, while $\mathbf{s}(\uu_k)$ takes any unit vector when $\|\uu_k\| = 0$. The friction magnitude function, $f$, is nonsmooth with respect to $\uu_k$ since $f(\|\uu_k\|) = 1$ when $\|\uu_k\|> 0$, and $f(\|\uu_k\|) \in [0,1]$ when $\|\uu_k\| = 0$. \del{These} \add{This} nonsmoothness would severely slow and even break convergence of gradient-based optimization.

\begin{figure}[t]
    \centering
    \includegraphics[width=0.4\linewidth]{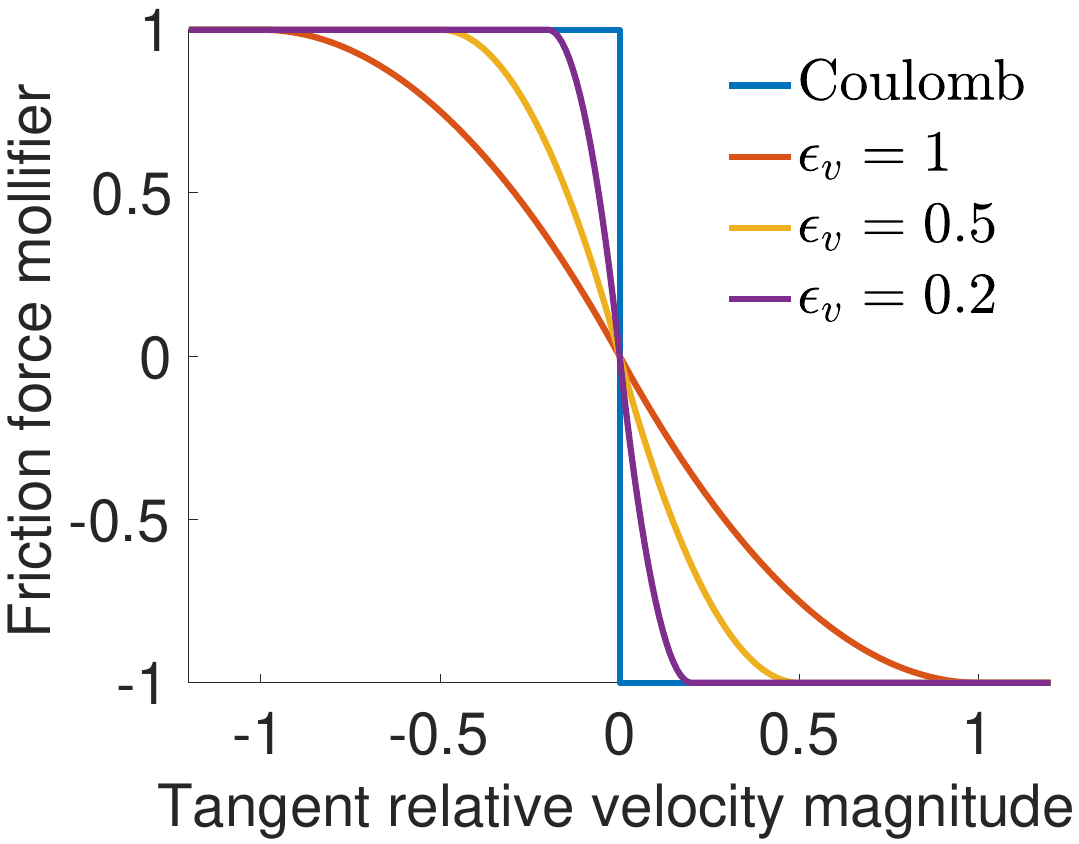}
    \caption{\textbf{Friction mollifier plotted with different $\epsilon_v$.} Decreasing $\epsilon_v$ asymptotically matches the discontinuous Coulomb friction model.}
    \label{fig:friction_plot}
\end{figure}

To enable efficient and stable optimization, the friction-velocity relation in the transition to static friction can be mollified by replacing $f$ with a smoothly approximated function. Following Li \etal\ \cite{li2020incremental}, we use
\begin{align}
f_1(y) =
\begin{cases}
 	-\frac{y^2}{\epsilon_v^2 \Delta t^2} + \frac{2y}{\epsilon_v \Delta t}, & y \in[0,\Delta t \epsilon_v)\\
 	1, & y \geq \Delta t \epsilon_v,
 \end{cases}
 \label{eq:fric_mollifier}
\end{align}
where $f_1'(\Delta t\epsilon_v) = 0$ and a velocity magnitude bound $\epsilon_v$ (in units of $m/s$) below which sliding velocities $\vv_k$ are treated as static is defined for bounded approximation error (Figure\ \ref{fig:friction_plot}).

Note that the velocity used in our friction model on the MPM side is the interpolated grid velocity at particle quadrature locations, rather than the particle velocity after grid-to-particle transfer. This makes the velocity seen by frictional forces independent from the choice of the particle-grid transfer scheme. This is important because for example in FLIP, the particle velocity does not reflect its displacement ($\vv_q^{n+1}\neq (\xx_p^{n+1}-\xx_p^n)/\Delta t$) and thus should not be used to define friction in an implicit solve.

However, challenges remain on incorporating friction into the optimization time integration. A major problem is that friction is not a conservative force and there is no well-defined potential such that taking the opposite of its gradient produces the frictional force. Therefore, following Li \etal\ \cite{li2020incremental}, we fix the friction constraint set $\mathcal{F}$ along with the normal force magnitude $\lambda$ and the tangent operator $T$ during the nonlinear optimization to the last updated value $\mathcal{F}^j = \mathcal{F}(x^j)$, $\lambda^j = \lambda(x^j)$, and $T^j=T(x^j)$, which then makes the lagged friction force integrable with the pseudo-potential
\begin{equation}
    D(x) = \sum_{k\in \mathcal{F}^j } \mu \lambda_k^j f_0(\|\bar{\uu}_k\|),
\end{equation}
where $\mathcal{F}^j$ is the set of all contact pairs with nonzero $\lambda_k^j$, $f'_0(y)=f_1(y)$, $\bar{\uu}_k = (T_k^j)^T (\Delta t v^n + \Delta t^2 ((1-\gamma) a^{n} + \gamma a^{n+1}))$ and so we have $-\nabla D(x) = -\sum_{k\in\mathcal{F}^j}\mu \lambda^j_k T^j_k f_1(\|\bar{\uu}_k\|) \> \mathbf{s}(\bar{\uu}_k)$, which is a semi-implicit discretization of the frictional force with lagged variables $\lambda^j_k$ and $T^j_k$.
Then we can iteratively alternate between the nonlinear optimization with fixed $\mathcal{F}$, $\lambda$, and $T$ given as
\begin{equation}
    \min_{x}: E(x) = \frac{1}{2} \|x - \hat{x}^{n}\|^2_M + 2 \alpha \beta \Delta t^2 \left( \Psi(x) + B(x) + D(x) \right),
    \label{eq:discretized_IP_withBandD}
\end{equation}
and friction update until convergence (Algorithm\ \ref{alg:main}). Although the friction convergence is not guaranteed for arbitrarily large time step sizes due to the nonlinearity and asymmetry of the problem, we have confirmed that all our experiments converge with the practical time step sizes applied (Section\ \ref{sec:experiments}).

\subsection{Irregular Boundaries for MPM}

In the BFEMP framework, an experimental setup with a subset or all of the FEM nodes prescribed with Dirichlet boundary conditions on their displacements can be applied to model irregular boundaries for MPM. This can not only resolve detailed boundary geometries even when the MPM grid is relatively \del{coarser} \add{coarse} (Section\ \ref{sec:fem_as_bc}), but also provide accurate and controllable friction on the boundary (Section\ \ref{sec:exp_rectangle_slope}).
\section{Nonlinear Optimization}\label{sec:optimization}

The time integration framework of BFEMP for one time step is outlined in Algorithm\ \ref{alg:main}. MPM particle-grid transfers are performed in the beginning (line 2) and the end (line 12). On the MPM grid and the FEM mesh, the minimization of incremental potential with lagged friction (line 7) is alternated with \add{the} friction update (line 9) until convergence to the fully implicit friction solution.

\begin{algorithm}[ht!]
\caption{BFEMP Time Integration}
\label{alg:main}
\begin{algorithmic}[1]
\Procedure{TimeIntegration}{$x_\text{F}^n$, $v_\text{F}^n$, $M_\text{F}$, $x_\text{P}^n$, $v_\text{P}^n$, $M_\text{P}$, $\Delta t$} \Comment{subscript P is for stacked particle variables}
    \State $x_\text{M}^n$, $v_\text{M}^n$, $M_\text{M}^n \leftarrow$ particleToGrid($x_\text{P}^n$, $v_\text{P}^n$, $M_\text{P}$) \Comment{Table\ \ref{tb:apic},\ \ref{tb:pic}, and\ \ref{tb:flip}}
    \State $\tilde{x}^{n+1}_\text{F} \leftarrow {x}^{n}_\text{F}$, $\tilde{x}^{n+1}_\text{M} \leftarrow {x}^{n}_\text{M}$ \Comment{for initial guess}
    \State $j \leftarrow 0$
    \State $\mathcal{F}^j$, $\lambda^j$, $T^j$\del{,} $\leftarrow$ computeFrictionOperator($\tilde{x}_\text{F}^{n+1}$, $\tilde{x}_\text{M}^{n+1}$) \Comment{Section\ \ref{sec:friction}}
    \Do
        \State \begin{math}\begin{bmatrix}\tilde{x}^{n+1}_\text{F} \\ \tilde{x}^{n+1}_\text{M}\end{bmatrix}\end{math},
        \begin{math}\begin{bmatrix}\tilde{v}^{n+1}_\text{F} \\ \tilde{v}^{n+1}_\text{M}\end{bmatrix}\end{math} 
        $\leftarrow$ 
        MinimizeIP(
            \begin{math}\begin{bmatrix}x^n_\text{F} \\ x^n_\text{M}\end{bmatrix}\end{math},
            \begin{math}\begin{bmatrix}v^n_\text{F}\\ v^n_\text{M}\end{bmatrix}\end{math},
            \begin{math}\begin{bmatrix}M_\text{F} & \\ & M^n_\text{M}\end{bmatrix}\end{math},
            $\Delta t$, $\mathcal{F}^j$, $\lambda^j$, $T^j$,
            \begin{math}\begin{bmatrix}\tilde{x}^{n+1}_\text{F} \\ \tilde{x}^{n+1}_\text{M}\end{bmatrix}\end{math}) \Comment{Algorithm\ \ref{alg:newtonOTI}}
        \State $j \leftarrow j+1$
        \State $\mathcal{F}^j$, $\lambda^j$, $T^j$\del{,} $\leftarrow$ computeFrictionOperator($\tilde{x}_\text{F}^{n+1}$, $\tilde{x}_\text{M}^{n+1}$) \Comment{Section\ \ref{sec:friction}}
    \doWhile{friction not converged} \Comment{Section\ \ref{sec:optimization}}
    \State $x^{n+1}_\text{F} \leftarrow \tilde{x}^{n+1}_\text{F}$, $v^{n+1}_\text{F} \leftarrow \tilde{v}^{n+1}_\text{F}$
    \State $x^{n+1}_\text{P}$, $v^{n+1}_\text{P}$ $\leftarrow$ gridToParticle($\tilde{x}^{n+1}_\text{M}$, $\tilde{v}^{n+1}_\text{M}$) \Comment{Table\ \ref{tb:apic},\ \ref{tb:pic}, and\ \ref{tb:flip}}
    \State \textbf{return} $x_\text{F}^{n+1}$, $v_\text{F}^{n+1}$, $x_\text{P}^{n+1}$, $v_\text{P}^{n+1}$
\EndProcedure
\end{algorithmic}
\end{algorithm}

\begin{algorithm}[ht!]
\caption{Line Search Method for Incremental Potential Minimization}
\label{alg:newtonOTI}
\begin{algorithmic}[1]
\Procedure{MinimizeIP}{$x^n$, $v^n$, $M^n$, $\Delta t$, $\mathcal{F}^j$, $\lambda^j$, $T^j$, $\bar{x}$}
    \State $x \leftarrow \bar{x}$ \Comment{initial guess}
    \State $E_\text{prev} \leftarrow E(x)$, $x_\text{prev} \leftarrow x$ \Comment{$E(x)$ defined in (\ref{eq:discretized_IP_withBandD}) also depends on $x^n$, $v^n$, $M^n$, $\Delta t$, $\mathcal{F}^j$, $\lambda^j$, $T^j$}
    \Do
        \State $H \leftarrow \text{computeProxyMatrix}(x)$ \Comment{applying projected Newton\ \cite{teran2005robust}}
        \State $p \leftarrow -H^{-1} \nabla E(x)$ \Comment{solved using CHOLMOD\ \cite{chen2008algorithm}}
        \State $\tau \leftarrow \text{initStepSize}(x)$ \Comment{line search filtering\ \cite{wachter2006implementation}}
        \Do \Comment{Armijo line search\ \cite{nocedal2006numerical}}
            \State $x \leftarrow x_\text{prev} + \tau p$
            \State $\tau \leftarrow \tau / 2$
        \doWhile{$E(x) > E_\text{prev}$}
        \State $E_\text{prev} \leftarrow E(x)$, $x_\text{prev} \leftarrow x$
    \doWhile{$\|p\|_\infty/\Delta t > \epsilon_d$} \Comment{Section\ \ref{sec:optimization}}
    \State $\tilde{x}^{n+1} \leftarrow x$, $\tilde{v}^{n+1} \leftarrow v^n + \frac{1}{\Delta t} (M^n)^{-1} ((\gamma - 1)\nabla E(x^n) - \gamma \nabla E(x))$ \Comment{Section\ \ref{sec:discretization}}
    \State \textbf{return} $\tilde{x}^{n+1}$, $\tilde{v}^{n+1}$
\EndProcedure
\end{algorithmic}
\end{algorithm}

Applying the projected Newton's method\ \cite{teran2005robust} for incremental potential minimization (Algorithm\ \ref{alg:newtonOTI}), we compute the proxy matrix $H$ by projecting the local Hessian of every elasticity, barrier, and friction stencil to its closest positive semi-definite form by zeroing out the negative eigenvalues, and then summing them up together with the mass matrix $M^n$ (line 5). The search direction $p$ is computed by factorizing $H$ and back-solving it on $-\nabla E(x)$ using CHOLMOD\ \cite{chen2008algorithm} (line 6). To obtain global convergence, the backtracking line search that ensures the decrease of energy is applied (line 8 to 11), starting from a large feasible step size that avoids interpenetration and deformation gradient degeneracy (line 7). After converging to a local optimum, velocity is updated with the newly obtained position (line 14) and returned together (line 15).
Here we use the infinity norm of the Newton increment (search direction $p$) in the unit of velocity ($m/s$) for the stopping criteria, which provides a 2nd-order approximation on the distance to the true solution. Similarly, friction convergence in Algorithm\ \ref{alg:main} is also determined this way, but with $\mathcal{F}$, $\lambda$, and $T$ computed using the current $x$.

Along Newton's search direction $p$, we compute the largest step size that will first result in a $0$ distance on any contact pair or a $0$ determinant on any deformation gradient. We then set the initial line search step size to be $0.9\times$ of this critical value.
The critical value for $0$ distance is computed via continuous collision detection (CCD)\ \cite{Li2021CIPC}, and for $0$ determinant it is just the smallest positive real root of a polynomial equation\ \cite{li2021lagrangian}.
This ensures that interpenetration or deformation gradient degeneracy could never happen throughout the simulation since the following backtracks always result in step sizes smaller than the critical value.

The numerical parameters in BFEMP all have physical meanings and directly control the extent of approximation to the continuous problem. 
To summarize, we have $\hat{d}$ (contact activation distance in $m$), $\epsilon_v$ (stick-slip velocity threshold in $m/s$), $\epsilon_d$ (Newton tolerance in $m/s$), and physical parameter $\kappa$ (barrier stiffness in $Pa$). \add{Here $\kappa$ also affects the convergence speed of the projected Newton method (Algorithm\ \ref{alg:newtonOTI}), but the convergence is always guaranteed eventually. In our experiments, we observed that setting $\kappa$ several orders of magnitude smaller than the average elasticity stiffness of the objects in the simulation can provide efficient convergence.}
\section{Numerical Simulations}\label{sec:experiments}

In this section, we provide 6 examples in 2D and 1 example in 3D to \add{verify} \del{validate} the contact model and the friction model in the proposed BFEMP approach. The numerical parameters $\hat{d}$, $\epsilon_v$ , $\epsilon_d$, and physical parameter $\kappa$ are all reported respectively in each experiment. If not mentioned otherwise, all elasticities are with the neo-Hookean model, and all particle-grid transfer schemes are with APIC. The visualized stresses are all the von Mises stress.

\subsection{Momentum and Energy Study}

\begin{figure}[t]
    \centering
    \includegraphics[width=\textwidth]{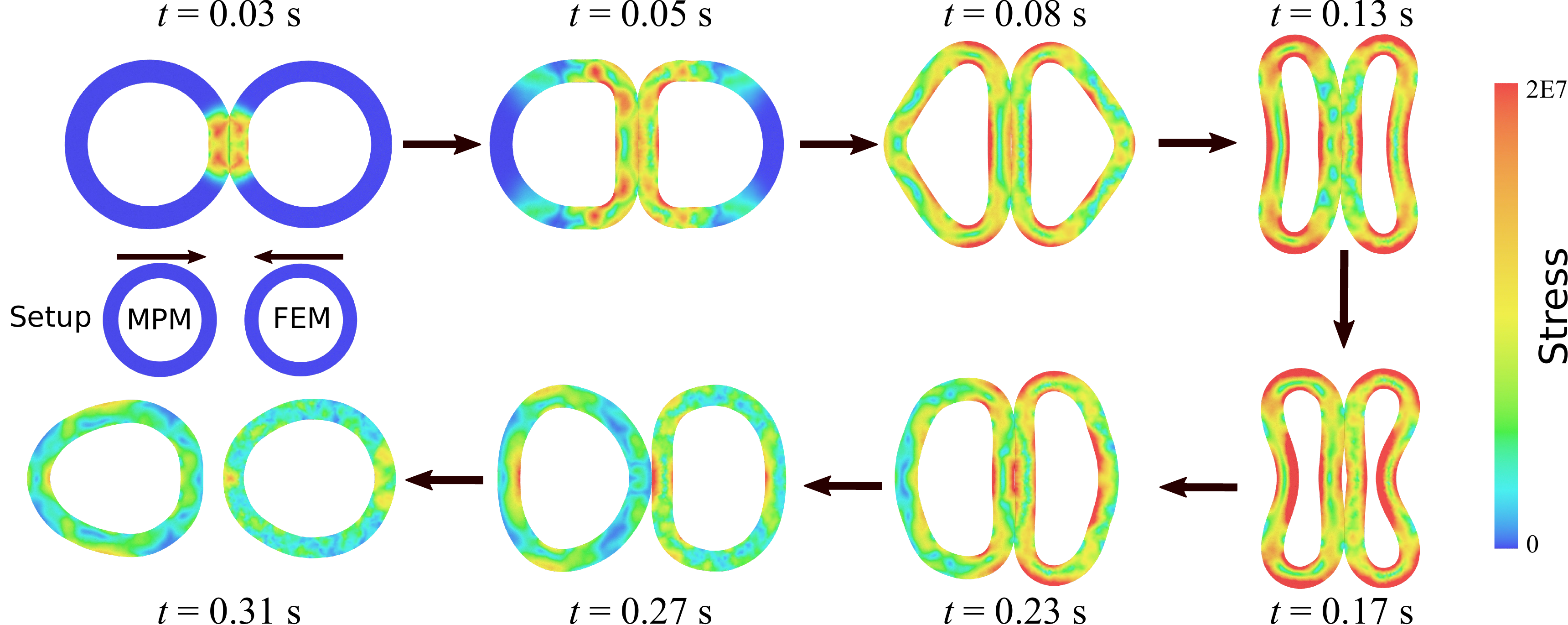}
    \caption{\textbf{Colliding rings.} The experiment setup and stress wave propagation over time.}
    \label{fig:elastic-rings-collision-sequence}
\end{figure}

\begin{figure}[t!]
    \centering
    \includegraphics[width=\textwidth]{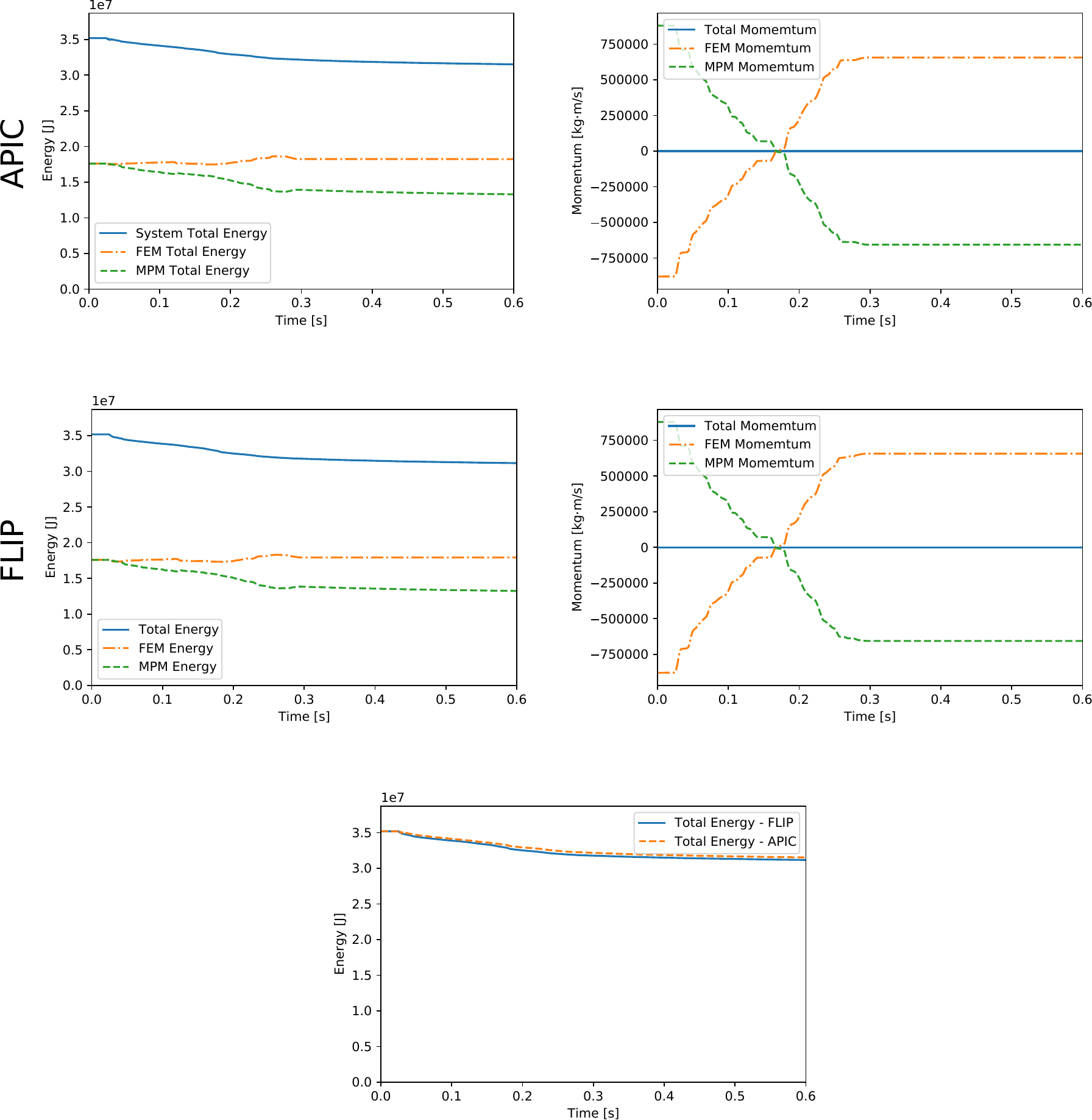}
    \caption{\textbf{Momentum and Energy Behavior.} The energy and momentum plot for APIC and FLIP transfer schemes.}
    \label{fig:elastic-rings-collision-plot}
\end{figure}

The collision between two elastic rings is simulated to \add{verify} \del{validate} the momentum and energy behavior of BFEMP and to demonstrate the robustness of our framework in handling large deformation. This example is modified from the MPM-MPM contact version in \cite{huang2010contact}.

The experiment setup is shown in the left-middle subfigure of Figure \ref{fig:elastic-rings-collision-sequence}. The two rings are identical except that the left ring is discretized with MPM, and the right one is discretized with FEM. The inner radius of the ring shape is $3m$, and the outer radius of the ring shape is $4m$. The Young's modulus is $E = 10^8 Pa$, the Poisson's ratio is $\nu = 0.2$, and the density is $\rho = 1000 kg/m^2$. The MPM ring is discretized by 20098 particles, where the grid spacing is $0.1m$. The FEM ring is discretized by 1830 vertices and 3310 triangles. The gravitational force and the frictional forces are not included. The two rings are placed $2m$ apart and then move towards each other with an initial speed of $40 m/s$. The contact active distance and the contact stiffness are set to $\hat{d} = 10^{-2} m$ and $\kappa = 10^7 Pa$ respectively. To minimize numerical dissipation, we use the Newmark time integrator with time step size $\Delta t = 2 \times 10^{-4} s$. The Newton tolerance is  $\epsilon_d = 10^{-6} m/s$.  We also compare the APIC and FLIP transfer schemes. Their differences in displacement and stress are small, so only the stress wave propagation with APIC is shown in Figure \ref{fig:elastic-rings-collision-sequence}. The energies and momenta over time are plotted for both APIC and FLIP in Figure \ref{fig:elastic-rings-collision-plot}.

The collision happens between $0.025s$ and $0.293s$. The symmetry of stress patterns is preserved during the collision. The system's total momentum is perfectly preserved with both transfer schemes. Part of the energy is lost during the collision: $8.57\%$ energy is lost with APIC, and $9.67\%$ with FLIP. After the rings are separated, the FEM ring preserves its energy over time, while the MPM ring gradually loses energy, primarily due to numerical dissipation in the particle-grid transfers.

\subsection{FEM as Contact Boundary for MPM}\label{sec:fem_as_bc}

\begin{figure}[t]
    \centering
    \includegraphics[width=\textwidth]{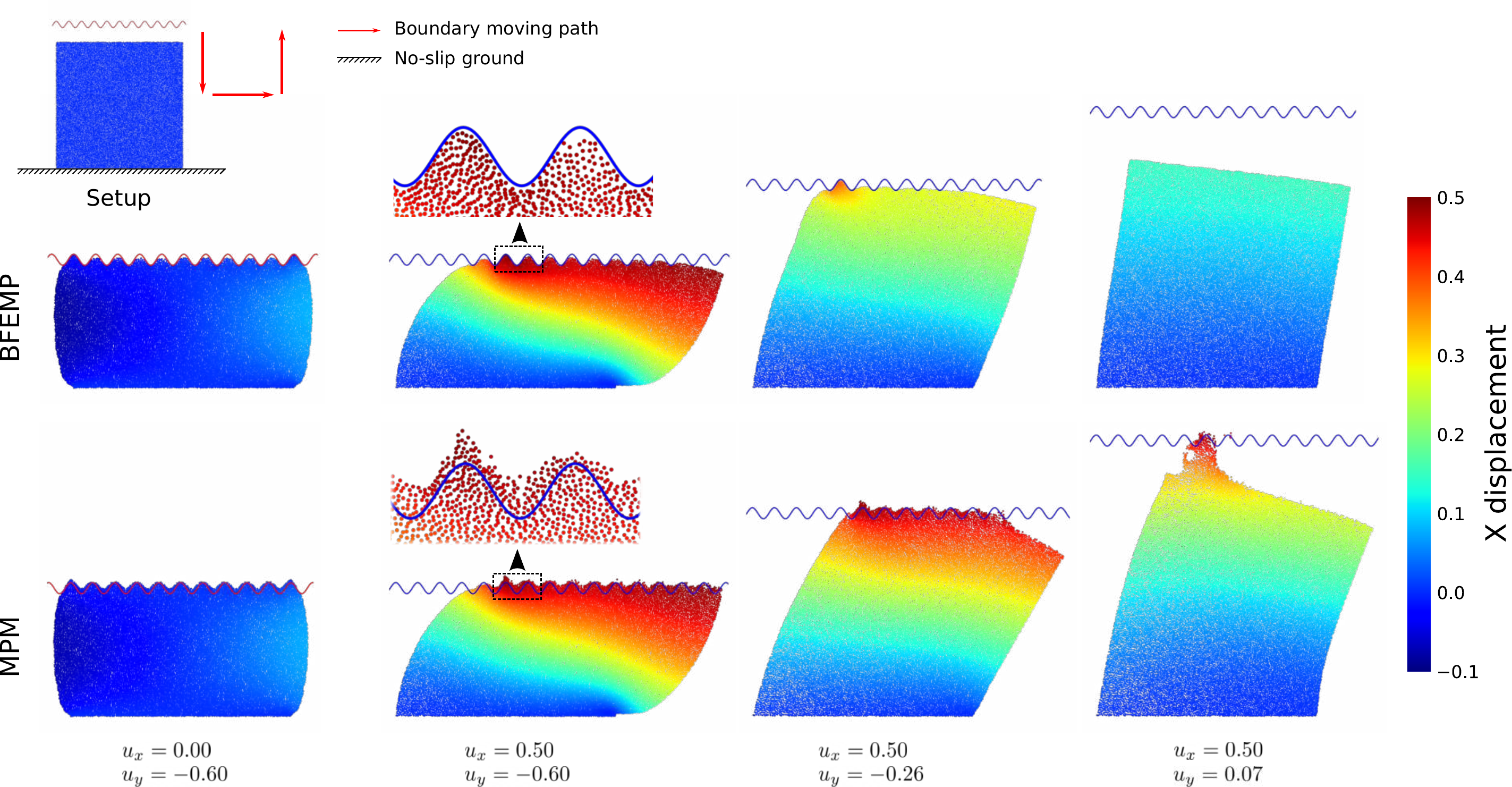}
    \caption{\textbf{FEM as Boundary Condition.} The friction-free interaction between a sine wave shape boundary and an MPM cube are simulated to compare BFEMP based slip boundary condition and traditional level set based slip boundary condition. $(u_x, u_y)$ is the displacement of the sine wave w.r.t. its initial position. BFEMP based slip boundary condition can guarantee non-penetration and doesn't have adhesive forces when it is separating from the object.}
    \label{fig:sine_wave}
\end{figure}

The guaranteed impenetrability between MPM particles and FEM boundaries makes BFEMP a natural strategy for enforcing kinematic separable boundary conditions in MPM simulations. Here we test the friction-free interaction between a sine wave shape boundary and an MPM cube. The BFEMP-based boundary condition is compared with a level-set based slip boundary condition, which enforces a zero normal relative velocity condition at each grid node inside the sine wave's level set, \add{i.e., at each time step, for those nodes that are within the level set, their normal velocities along the level set interface are prescribed, so that the original unconstrained optimization \ref{eq:discrete_IP} for the time integration are solved with these equality constraints. }

The experiment setup is shown in the left-top subfigure of Figure \ref{fig:sine_wave}. A $1m \times 1m$ elastic box with Young's modulus $E = 10^6 Pa$, Poisson's ratio $\nu = 0.2$ and $\rho = 1000 kg/m^2$ is placed on a no-slip ground. It is discretized by 21026 particles, with grid spacing $0.02m$. A sine wave boundary is placed $0.2m$ above the box, whose contour is determined by $y = \frac{1}{40} \cos \frac{2 \pi}{0.1} x $. For BFEMP, the sine wave boundary condition is discretized by a FEM mesh with prescribed displacements at each time step. While for MPM, it is described by an analytical level set. The sine wave boundary first moves $0.6m$ downwards, then $0.5m$ to the left, and finally upwards until separation. The moving speed is $1m/s$ all the way. The contact active distance and the contact stiffness is set to $\hat{d} = 10^{-3} m$ and $\kappa = 10^4 Pa$ respectively. The implicit Euler time integration with time step $h = 10^{-3} s$ is used. The Newton tolerance is set to $\epsilon_d = 10^{-4} m/s$.

As shown in Figure \ref{fig:sine_wave}, the BFEMP-based boundary condition more accurately resolves the complex boundary geometry without exhibiting any numerical adhesive forces when the boundary is separating from the cube. With the level-set based slip condition, particles will penetrate the boundary because the boundary condition is only defined on the MPM grid in a ``smeared out'' manner. The numerical adhesive force comes from that at each time step, the grid with slip condition is locked within some plane. On the contrary, with BFEMP, MPM particles can freely move around outside the FEM mesh.

\subsection{Brazilian Disk Test}

\begin{figure}[t]
    \centering
    \includegraphics[width=\textwidth]{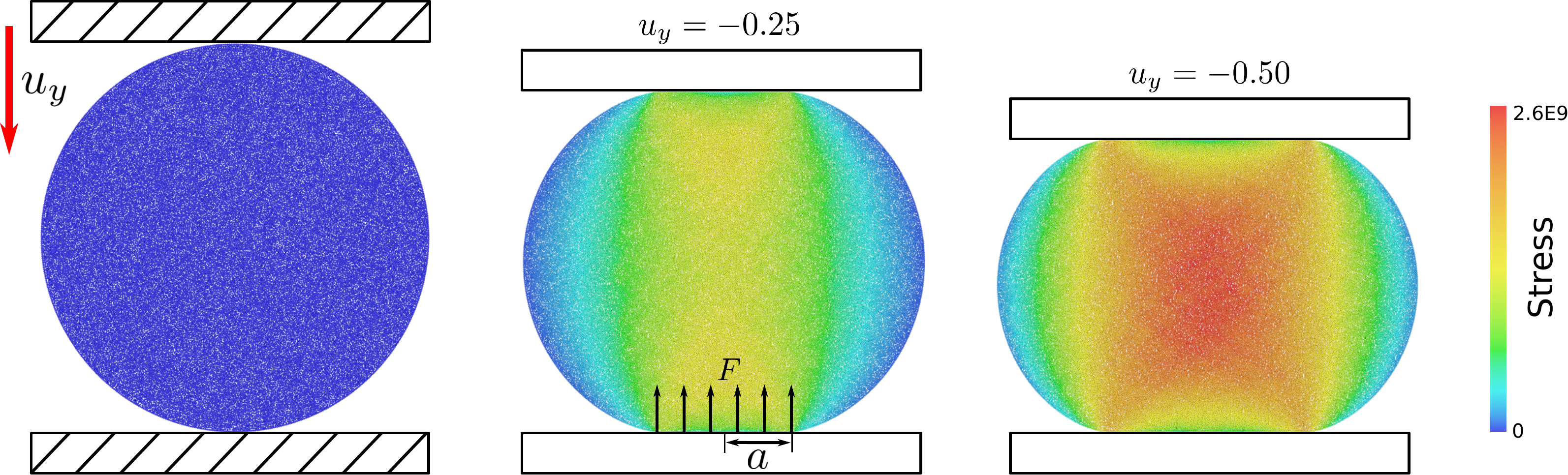}
    \caption{\textbf{Brazilian Disk Test.} The experiment setup and the compression procedure are shown here. The contact force and contact radius are illustrated in the middle figure.}
    \label{fig:contact_force_test}
\end{figure}

\begin{figure}[ht!]
    \centering
    \includegraphics[width=\textwidth]{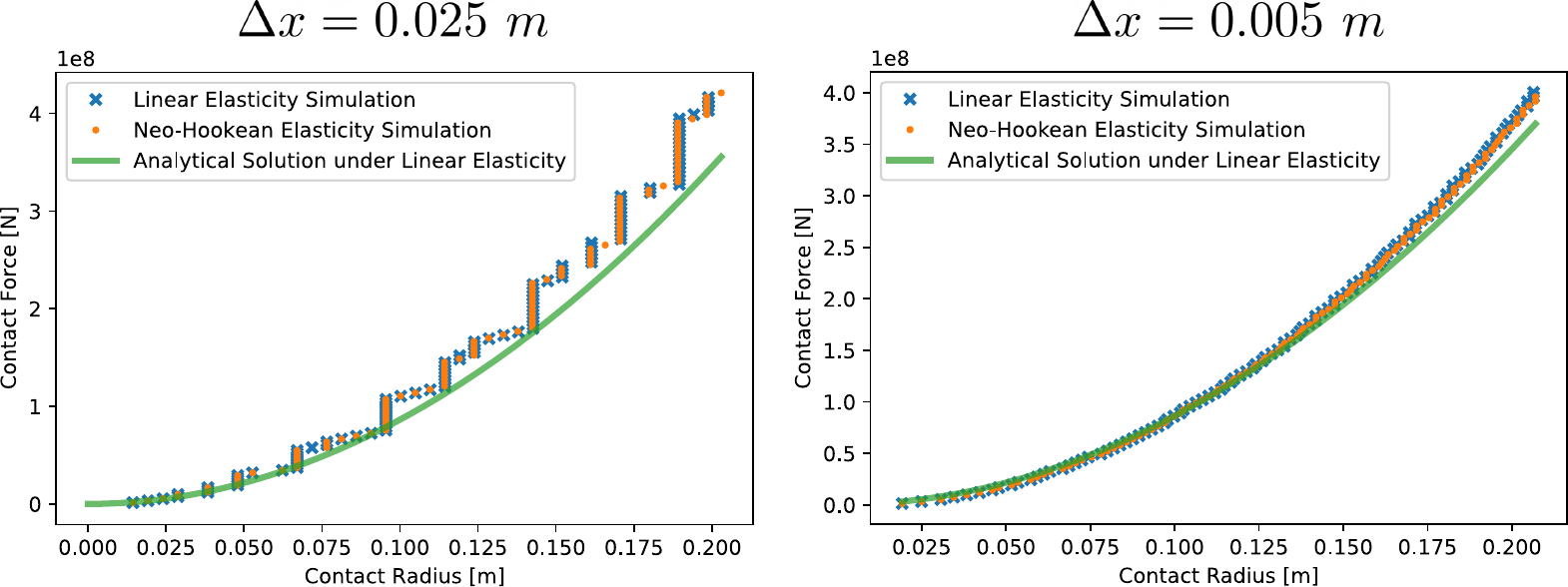}
    \caption{\textbf{Brazilian Disk Test.} Within the small deformation range, our contact model fits well with Hertzian contact theory. \add{The non-smoothness of the measured radius from the simulation results can be alleviated as the resolution increases.}}
    \label{fig:contact_force_test_plot}
\end{figure}

To \add{verify} \del{validate} the accuracy of the contact model, BFEMP is studied on the Brazilian disk test, which is a special case of the plane Hertzian contact problem \cite{barber2009contact, liu2020ilsmpm}. The Brazilian disk test can be used for tensile strength testing, which involves a 2D elastic disk squeezed between two rigid objects. We use a fixed rigid plate and a moving rigid plate to simulate the compression procedure. According to the Hertzian contact model, the contact force $F$ and the contact radius $a$ have the following relation:

\begin{equation}
    F = \frac{\pi}{4} \frac{E}{1 - \nu^2} \frac{a^2}{R}.
\end{equation}
The contact force and the contact radius are illustrated in Figure \ref{fig:contact_force_test_plot}.

In this experiment, the radius of the MPM disk is $1m$. It is composed of 42920 particles with MPM grid spacing $\Delta x = 0.025m$. The Young's modulus is $E = 10^{10}Pa$, and the Poisson's ratio is $\nu= 0.3$. To reduce the \add{inertial} effect, \del{of the inertial force,} \add {we artificially decrease the density of the material, which } \del{the density} is set to $\rho = 100 kg/m^2$\del{, which is much smaller than its Young's modulus}. The contact active distance is $\hat{d} = 10^{-4}m$ and the contact stiffness is $\kappa = 10^4 Pa$. The two plates are discretized with FEM. Each of them is composed of four vertices and two triangles. The fixed plate is placed $\hat{d}$ below the disk, and the moving plate is placed $\hat{d}$ above the disk. The constant velocity $0.1m/s$ of the moving plate is enforced by prescribing its displacements at each time step. The simulation is performed with implicit Euler time integration with time step size $h = 10^{-2}s$ and the Newton tolerance is $\epsilon_d = 10^{-8} m/s$. Friction coefficient $\mu = 1$ is used to prevent the disk from slipping.

The Hertzian model requires to measure the contact radius $a$. Following \cite{liu2020ilsmpm}, we use half of the horizontal range of the particles within the contact distance around the bottom FEM plate to approximate it. Here we test both linear elasticity and neo-Hookean elasticity. The compression procedure in Figure \ref{fig:contact_force_test} is visualized for the linear elasticity case. The $(a, F)$ data points within the small deformation range from the two simulations and the analytical $F-a$ relation from the Hertzian contact model are plotted in Figure \ref{fig:contact_force_test_plot}. The non-smoothness of the \add{measured radius from the simulation results} \del{simulated data} is due to the inaccurate approximation of \add{the radius} $a$ through a finite number of particles. \add{This non-smoothness can be alleviated as the resolution increases.} Despite that, we observe a qualitative match between the simulated data and the Hertzian contact theory.

\subsection{Critical Value of Friction Coefficient} \label{sec:exp_rectangle_slope}

\begin{figure}[t]
    \centering
    \includegraphics[width=\linewidth]{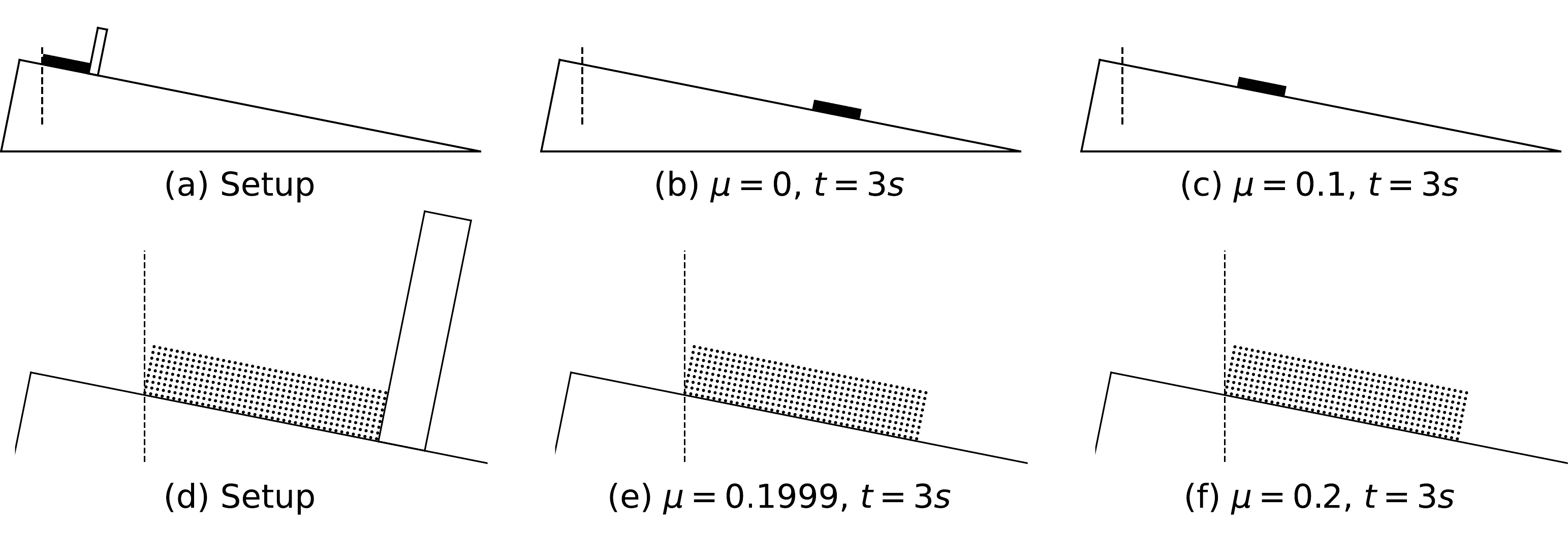}
    \caption{\textbf{Critical Value of Friction Coefficient.} (a,d) initial configuration with the extra support; (b,c,e) results at $t=3s$ of $\mu=0$, $0.1$, and $0.1999$, sliding distances all matching analytical solutions; (f) the result at $t=3s$ of $\mu=0.2$, static solution with sliding error bounded by $\epsilon_v$.}
    \label{fig:slope_test}
\end{figure}

\begin{figure}[ht]
    \centering
    \includegraphics[width=\linewidth]{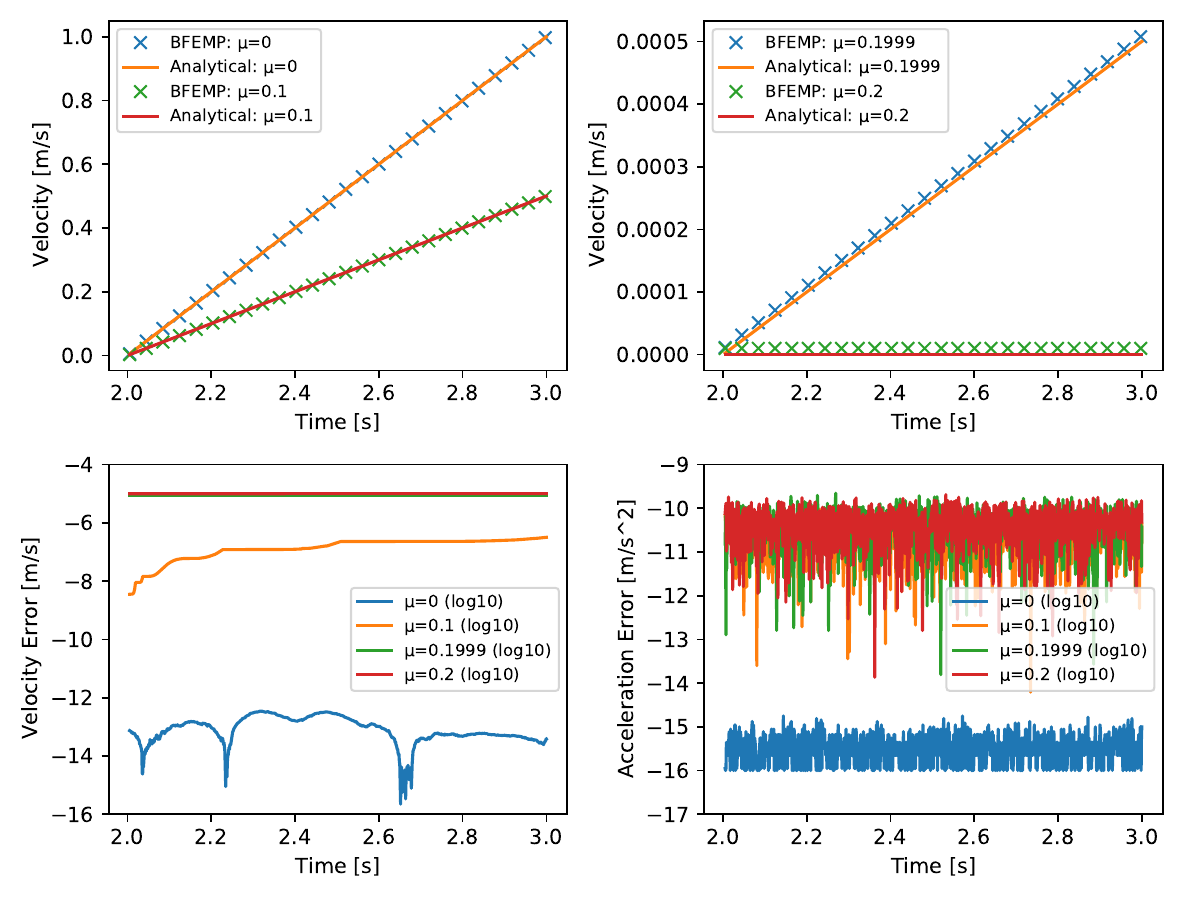}
    \caption{\textbf{Critical Value of Friction Coefficient.} At all friction coefficients, including $\mu=0$ (no friction), $\mu=0.1$, $\mu=0.1999$ ($99.95\%$ of the critical value), and $\mu=0.2$ (the critical value), the velocities and accelerations over the releasing period ($2s$ to $3s$) are all accurately matching the analytical solutions.}
    \label{fig:slope_test_plot}
\end{figure}

To \add{verify} \del{validate} the accuracy of BFEMP's friction model, an experiment with a stiff MPM box resting or sliding on a fixed FEM slope (or BFEMP's friction-controllable boundary condition) with a certain friction coefficient is created. When a rigid box is placed on a slope with zero initial velocity, its acceleration has the following analytical form:
\begin{equation}
a_x = 
\begin{cases}
    g (\sin\theta - \mu \cos\theta), & \theta \ge \tan \theta,\\
    0, & \theta < \tan \theta,
\end{cases}
\end{equation}
where $\mu$ is the friction coefficient between the box and the slope, $g$ is the gravity acceleration, $\theta \in [0, \pi/4)$ is the inclined angle of the slope.
Experiments show that BFEMP's friction model matches analytical solutions on sliding dynamics and critical value of friction coefficient both with bounded and small approximation error.

The initial configuration of this example is obtained by placing the MPM box $\hat{d}$ away from the slope, placing another fixed plane perpendicular to the slope on the side of the box where it may slide (also $\hat{d}$ away), and then simulate under gravity ($g=5.10m/s^2$) without friction until the box becomes static (Figure\ \ref{fig:slope_test}a). After obtaining the initial configuration, the slope test simulation is performed without the extra plane and with multiple different friction coefficients for each test (Figure\ \ref{fig:slope_test}b,c,d).

Here the MPM box is $0.1m\times 0.02m$, composed of 369 particles (grid $dx=0.005$) with density $\rho=100kg/m^2$, Young's modulus $E=4.0\times10^{12}Pa$ and Poisson's ratio $\nu=0.2$. Slopes with friction coefficient $\mu=0$, $0.1$, $0.1999$, and $0.2$ have been tested, all with contact active distance $\hat{d}=0.001m$, contact stiffness $\kappa=10^6Pa$, static friction velocity threshold $\epsilon_v=10^{-5}m/s$, and with the lagged normal forces in friction iteratively updated until converging to a solution with fully-implicit friction. All simulations are using implicit Euler time integration with time step size $h=0.001s$, and the Newton tolerance is set to $\epsilon_d=10^{-8}m/s$.

With sliding velocity and acceleration of the box's center of mass plotted over time (Figure\ \ref{fig:slope_test_plot}), they have all been shown to well match analytical solutions within $0.01\%$ relative errors. Even for $\mu=0.1999$ ($99.95\%$ that of the critical coefficient), the sliding behavior can still be accurately captured. For $\mu=0.2$, it is also confirmed that the acceleration vanishes, and the velocity throughout the simulation is around $\epsilon_v$, the static friction velocity threshold in BFEMP's approximation to provide the static friction force in the same magnitude as dynamic friction.

\subsection{Convergence under Refinement}

\begin{figure}[t]
    \centering
    \includegraphics[width=0.8\textwidth]{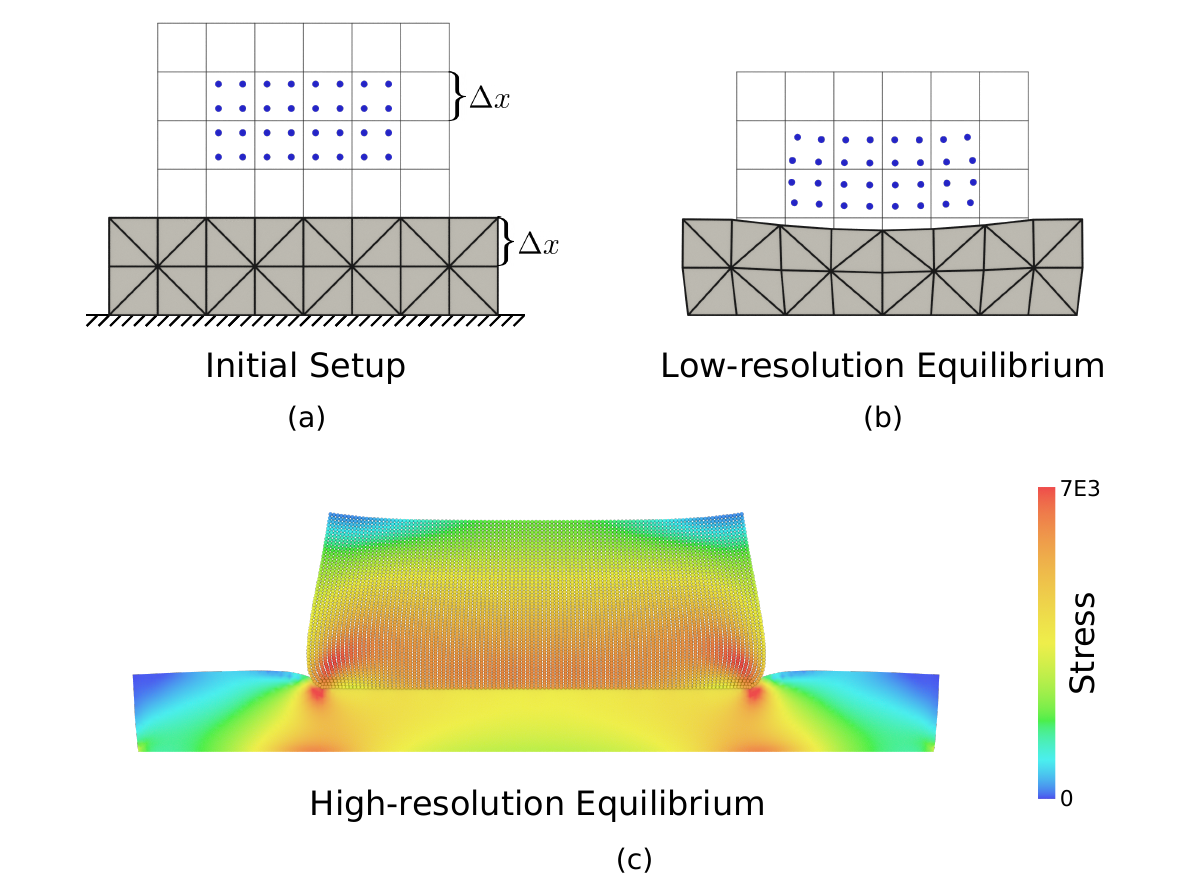}
    \caption{\textbf{Convergence under Refinement.} (a) Experiment setup. (b) The final equilibrium under low resolution. (c) The final equilibrium under high resolution. Stress pattern is visualized.}
    \label{fig:convergence_study_illustration}
\end{figure}

\begin{figure}
    \centering
    \includegraphics[width=0.6\textwidth]{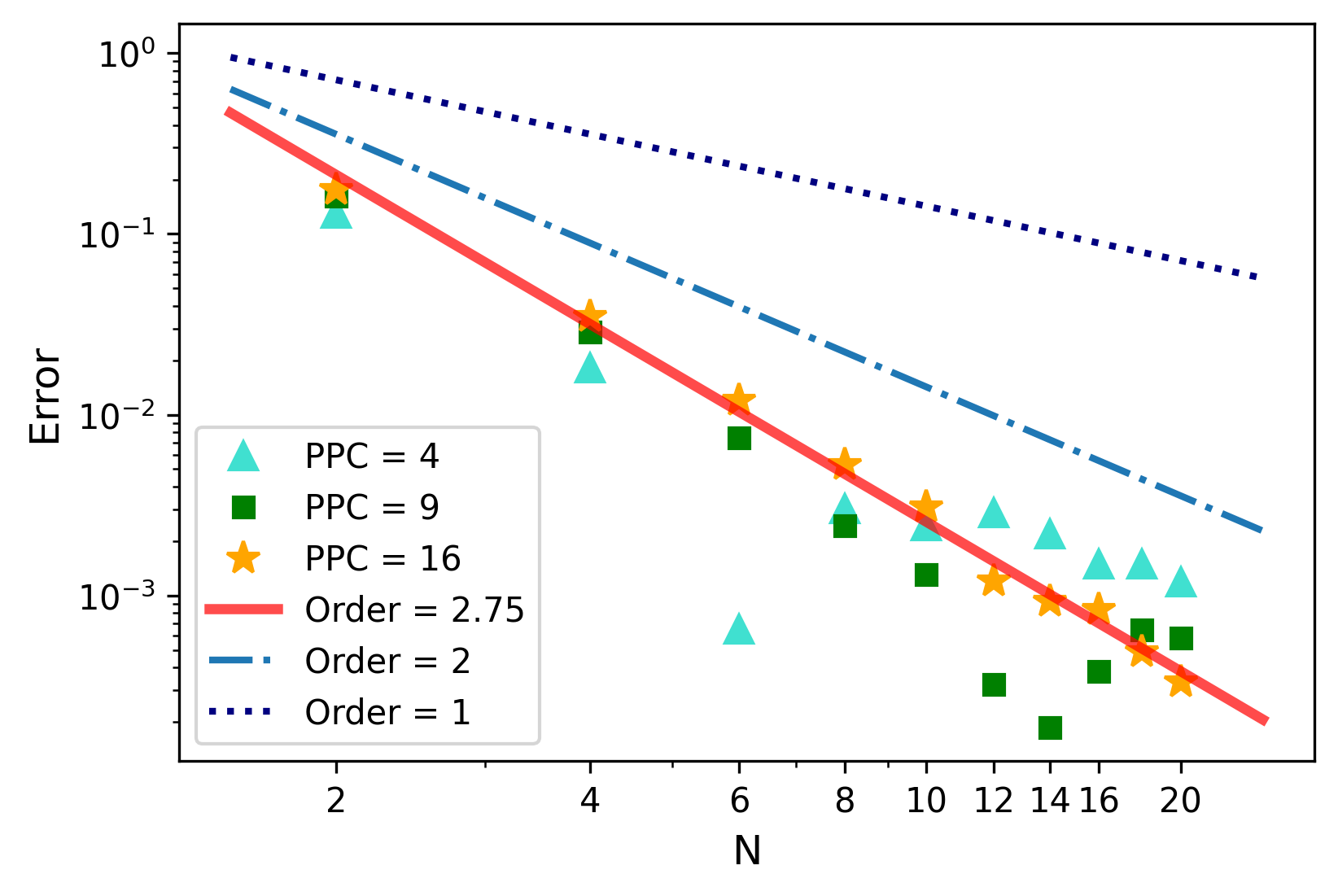}
    \caption{\textbf{Convergence under Refinement.} Higher PPC can reduce the noise in the convergence curve at higher resolutions. BFEMP with PPC = 16 achieves a convergence order of 2.75 to high-resolution result. Convergence curves with order 1 and order 2 are also plotted for reference.}
    \label{fig:convergence_study_plot}
\end{figure}

To \add{verify} \del{validate} the convergence under refinement property of BFEMP, an example with a soft MPM box stacking on a soft FEM box is created. A series of experiments with increasing resolutions and decreasing contact active distances are simulated to study the convergence rate under refinement. Results show that BFEMP can achieve a second-order convergence rate.

The initial configuration is illustrated in Figure \ref{fig:convergence_study_illustration} (a). The MPM box is with size $2m \times 1m$, Young's modulus $E = 4 \times 10^4 Pa$, Poisson's ratio $\nu = 0.4$ and density $\rho = 10^3 kg/m^3$. The particles are sampled regularly within each cell by placing each particle on the center of a sub-cell. The FEM box below is with size $4m \times 1m$, Young's modulus $E = 4 \times 10^4 Pa$, Poisson's ratio $\nu = 0.4$ and density $\rho = 10^2 kg/m^2$. The minimal edge length of the FEM mesh and the grid spacing of MPM are with the same value ($\Delta x$) in each experiment. The MPM box initially is placed $\Delta x$ above the FEM box and then simulated under gravity ($g = 10m/s^2$) until no oscillation is observed. To accelerate simulation to reach its final static state, PIC transfer scheme and implicit Euler time integration with large time step sizes (up to CFL limit for MPM) are used. The Newton tolerance is set to $\epsilon_d = 10^{-9} m/s$. The contact stiffness is set to $\kappa = 10^6 Pa$ for all experiments. Figure \ref{fig:convergence_study_illustration} (b) and Figure \ref{fig:convergence_study_illustration} (c) show the final equilibria under low resolution and high resolution respectively.

To examine the convergence rate of displacement to high-resolution results, the example is refined with $\Delta x = \frac{1}{N}$ and $\hat{d} = \frac{1}{N^2}$, where $N$ iterates all positive even numbers smaller than or equal to $20$. The reference high-resolution result is choose as with $N = 30$. The error is defined as the difference in height of the center of mass of the whole domain (with both FEM and MPM domains) between each testing resolution and the high-resolution reference. Due to quadrature error in MPM \cite{de2019material}, we also experiment with three different particle per cell (PPC) values: 4, 9, and 16. The three error sequences are plotted in Figure \ref{fig:convergence_study_plot}. As observed from the plot, a higher PPC value can reduce the noise in the convergence curve. The error sequence with PPC 16 almost falls into line. Under this setting, BFEMP achieves a convergence order of $2.75$.

\subsection{Buckling Behaviours under Different Friction Coefficients}
This example tests frictions between two semi-circular rings with large deformation. The two semi-circular rings are stacked together. As the outer semi-circular ring is compressed, different buckling patterns of the inner semi-circular ring under different friction coefficients are observed. This example is modified from the version with FEM-FEM contact in \cite{zimmerman2018febio}.

\begin{figure}[t]
    \centering
    \includegraphics[width=\textwidth]{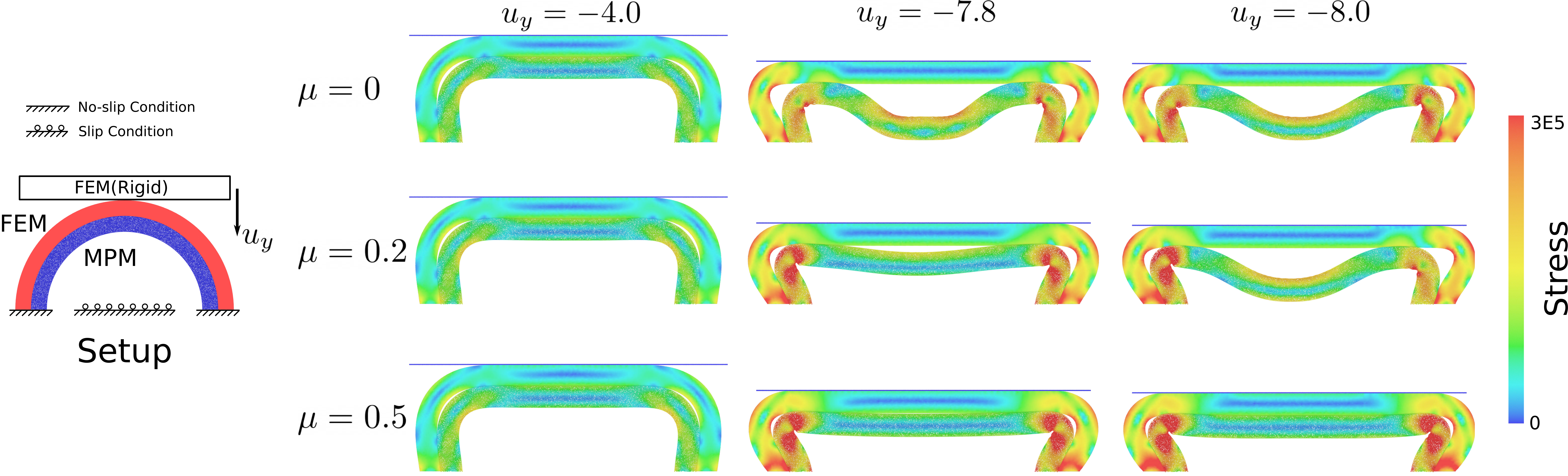}
    \caption{\textbf{Buckling Behaviours under Different Friction Coefficients.} The experiment setup is illustrated on the left. Under different friction coefficients, the buckling appears at different vertical displacements ($u_y$).}
    \label{fig:shell_friction}
\end{figure}

The experiment setup is shown in Figure \ref{fig:shell_friction}. The outer semi-circular ring with outer radius $14m$ and inner radius $12m$ is discretized by FEM with 2714 vertices and 5067 triangles. The inner semi-circular ring with outer radius $11.99m$ and inner radius $10m$ is discretized by MPM with 10261 particles with grid spacing $0.25m$. The two semi-circular rings are both with Young's modulus $E = 10^6 Pa$, Poisson's ratio $\nu = 0.3$ and density $\rho = 100 kg/m^2$. One FEM plate is placed $10^{-3}m$ above the outer semi-circular ring. The displacement of this plate is prescribed to follow a rigid linear motion with a constant downward velocity $1 m/s$. Large friction ($\mu = 10$) between the plate and the outer semi-circular ring is activated so that the plate can be viewed as a BFEMP based no-slip boundary condition. The feet of two semi-circular rings are fixed using the level set-based no-slip boundary condition. Another level set-based slip boundary condition is added at the bottom middle below the semi-circular rings to prevent the inner semi-circular rings from colliding into the ground. The contact active distance and the contact stiffness are set to $\hat{d} = 10^{-3}m$ and $\kappa = 10^5 Pa$. The static friction velocity threshold is set to $\epsilon_v = 10^{-5} m/s$. Implicit Euler time integrator with time step size $h = 10^{-2} s$ and Newton tolerance $\epsilon_d = 10^{-4} m/s$ are used.

We vary the friction coefficient between the two semi-circular rings from $\mu = 0$, $\mu = 0.2$ and $\mu = 0.5$. The compression procedure is visualized in Figure \ref{fig:shell_friction}. In the beginning, there is little difference between the three settings. As the FEM plate moves further down, the inner MPM semi ring under the friction-free setting is buckled first as expected. Friction with $\mu = 0.2$ lags the appearance of the buckling. For the large friction case with $\mu = 0.5$, no buckling happens at all.

\subsection{3D Twist with Friction}

\begin{figure}[ht!]
    \centering
    \includegraphics[width=0.8\textwidth]{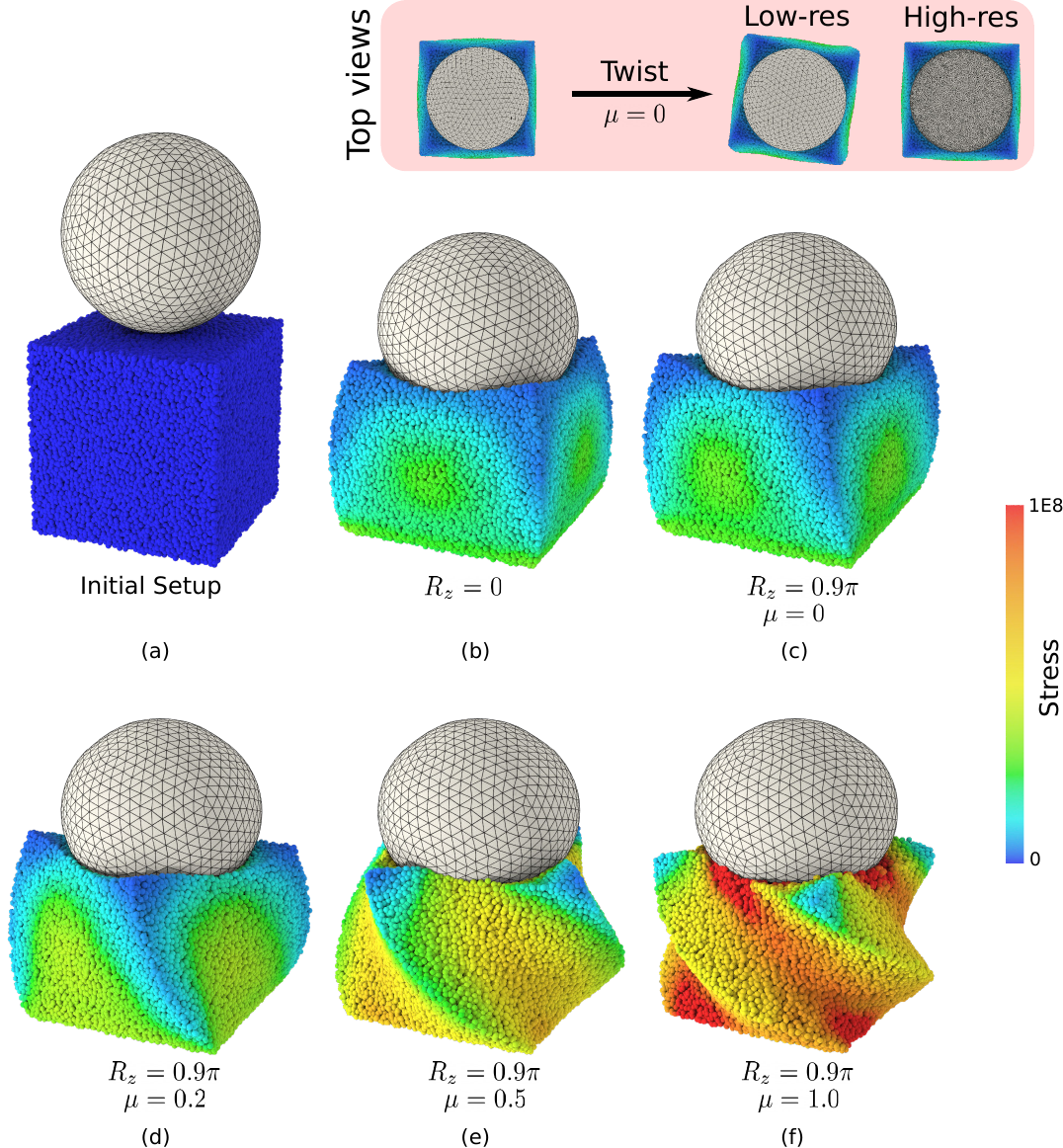}
    \caption{\textbf{3D Twist with Friction.} (a) Initial setup: The FEM spherical shell is placed $\hat{d}$ above the MPM cube; (b) Before twist procedure, the spherical shell is controlled to press down $0.5m$; (c, d, e, f) Equilibria under different friction coefficients when the shell stops rotating. \add{The nonzero rotation angle with $\mu = 0$ is caused by the non-smoothness of the contacting interfaces, which will decrease as the resolution increases (top views).}}
    \label{fig:3d_twist}
\end{figure}

\begin{figure}[ht!]
    \centering
    \includegraphics[width=0.6\textwidth]{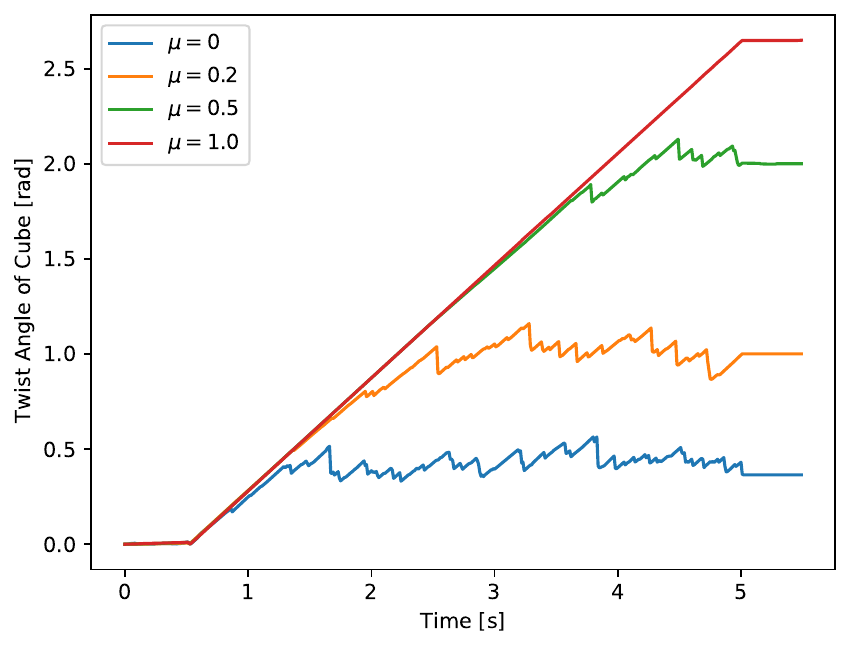}
    \caption{\textbf{3D Twist with Friction.} The average twist angle around the z-axis of the top center of the cube. The twist procedure happens between $0.5s$ and $5s$. The slipping between the cube and the shell appears at different time points under different friction coefficients.}
    \label{fig:3d_twist_angle_seq}
\end{figure}

To test BFEMP's contact and friction model in 3D, a twist test between a FEM spherical shell and an MPM cube is conducted. The FEM shell is controlled to exert a constant twist speed. Under different friction coefficients, it is expected to observe different maximal twist angles on the MPM cube. This example is modified from the version with FEM-FEM contact in \cite{zimmerman2018febio} as well.

The initial setup is illustrated in Figure \ref{fig:3d_twist} (a). The MPM cube with an edge length of $1m$ is placed $10^{-2} m$ below the FEM shell. It is discretized by 90929 particles, where the grid spacing is $0.0625m$. The Young's modulus is $10^8 Pa$. The Poisson's ratio is 0.4. And the density is $100 kg/m^3$. The FEM shell with inner radius $0.45m$ and outer radius $0.5m$ is discretized by 1364 points and 3920 tetrahedra. The Young's modulus is $10^{10} Pa$. The Poisson's ratio is 0.4. The density is $10^4 kg/m^3$. The displacements of the top part of the shell are prescribed to follow a rigid motion to exert downward compression and constant-speed twist: it first compresses down with a constant speed $0.5m/s$ for $1s$ (Figure \ref{fig:3d_twist} (b)) and then rotates around the z-axis with a constant angular velocity $\frac{\pi}{5}$ for $4.5s$. The contact active distance and the contact stiffness are set to $\hat{d} = 10^{-2}m$ and $\kappa = 10^7 Pa$. For settings with frictions, the static friction velocity threshold is set to $\epsilon_v = 10^{-3} m/s$. The simulation uses the implicit Euler time integrator with the time step size $h = 10^{-2} s$. The Newton tolerance is set to $\epsilon_d = 10^{-3} m/s$.

With different friction coefficients, the slipping between the shell bottom and the top center of the cube happens after different twist angles $R_z$. The final equilibria when the shell stops twisting are visualized in Figure \ref{fig:3d_twist} (c) (d) (e) (f). The twist angles of the top center part of the cube are plotted in Figure \ref{fig:3d_twist_angle_seq}. Since the pressure forces are between triangles and particles, the interface between the shell and the cube is not perfectly smooth. This roughness results in that the slipping happens when $R_z = 0.1\pi$ in the friction-free settings. \add{The final rotation angle should decrease as the resolution increases. To verify this, we increase the resolution of the FEM mesh and compare the final equilibria in the original setting and the higher-resolution setting. The top views are attached in Figure \ref{fig:3d_twist}, which shows that, with higher resolution, the final state of the cube is close to the initial state before the twisting.}  For $\mu = 0.2$ and $\mu =  0.5$, the slipping happens when the twists angles are around $R_z = 0.3 \pi$ and $R_z = 0.7 \pi$ respectively. With $\mu = 1.0$, there is no slipping between the shell and the cube.

\section{Conclusion}

In this paper, we proposed a new method for \add{monolithically coupling} \del{modeling the coupling through frictional contact between} an MPM domain and a FEM domain for elastodynamics \add{through frictional contact}. By approximating the non-interpenetration constraint with a barrier energy term and performing time integration using a variational formulation, our method guarantees that no particles will penetrate into the FEM mesh throughout the simulation. Furthermore, when the displacement of the FEM domain is fully prescribed, BFEMP reduces to an explicit mesh-based boundary treatment for MPM. Through numerical experiments validating the energy behavior, robustness, stability, and accuracy, we demonstrated the advantages of the proposed method.

\paragraph{Limitations and Future Works}
\add{For our current formulation, when MPM particles get very close to the FEM boundary, there is in fact a small portion of overlap between FEM and MPM domains even when there is no interpenetration. This is because MPM particles represent a region of the domain. Although the overlapping area vanishes under spatial refinement, it would still be interesting to also consider the size and deformation of the regions when defining the distance constraints. In addition, it would be meaningful to extend our framework to support cutting of MPM solids by FEM meshes, where the MPM particles on different sides of the FEM mesh should not communicate with each other even when the FEM mesh is much thinner than the MPM kernel.}

\bibliographystyle{elsarticle-num} 
\bibliography{references_arxiv}

\begin{thebibliography}{10}
\expandafter\ifx\csname url\endcsname\relax
  \def\url#1{\texttt{#1}}\fi
\expandafter\ifx\csname urlprefix\endcsname\relax\def\urlprefix{URL }\fi
\expandafter\ifx\csname href\endcsname\relax
  \def\href#1#2{#2} \def\path#1{#1}\fi

\bibitem{sulsky1995application}
D.~Sulsky, S.-J. Zhou, H.~L. Schreyer, Application of a particle-in-cell method
  to solid mechanics, Computer physics communications 87~(1-2) (1995) 236--252.

\bibitem{de2019material}
A.~de~Vaucorbeil, V.~P. Nguyen, S.~Sinaie, J.~Y. Wu, Material point method
  after 25 years: theory, implementation and applications, Submitted to
  Advances in Applied Mechanics (2019) 1.

\bibitem{harlow1962particle}
F.~H. Harlow, The particle-in-cell method for numerical solution of problems in
  fluid dynamics, Tech. rep., Los Alamos Scientific Lab., N. Mex. (1962).

\bibitem{brackbill1988flip}
J.~U. Brackbill, D.~B. Kothe, H.~M. Ruppel, Flip: a low-dissipation,
  particle-in-cell method for fluid flow, Computer Physics Communications
  48~(1) (1988) 25--38.

\bibitem{hughes2012finite}
T.~J. Hughes, The finite element method: linear static and dynamic finite
  element analysis, Courier Corporation, 2012.

\bibitem{de2020total}
A.~de~Vaucorbeil, V.~P. Nguyen, C.~R. Hutchinson, A total-lagrangian material
  point method for solid mechanics problems involving large deformations,
  Computer Methods in Applied Mechanics and Engineering 360 (2020) 112783.

\bibitem{de2021modelling}
A.~de~Vaucorbeil, V.~P. Nguyen, Modelling contacts with a total lagrangian
  material point method, Computer Methods in Applied Mechanics and Engineering
  373 (2021) 113503.

\bibitem{zhang2016material}
X.~Zhang, Z.~Chen, Y.~Liu, The material point method: a continuum-based
  particle method for extreme loading cases, Academic Press, 2016.

\bibitem{bardenhagen2004generalized}
S.~G. Bardenhagen, E.~M. Kober, The generalized interpolation material point
  method, Computer Modeling in Engineering and Sciences 5~(6) (2004) 477--496.

\bibitem{zhang2011material}
D.~Z. Zhang, X.~Ma, P.~T. Giguere, Material point method enhanced by modified
  gradient of shape function, Journal of Computational Physics 230~(16) (2011)
  6379--6398.

\bibitem{sadeghirad2011convected}
A.~Sadeghirad, R.~M. Brannon, J.~Burghardt, A convected particle domain
  interpolation technique to extend applicability of the material point method
  for problems involving massive deformations, International Journal for
  numerical methods in Engineering 86~(12) (2011) 1435--1456.

\bibitem{gan2018enhancement}
Y.~Gan, Z.~Sun, Z.~Chen, X.~Zhang, Y.~Liu, Enhancement of the material point
  method using b-spline basis functions, International Journal for numerical
  methods in engineering 113~(3) (2018) 411--431.

\bibitem{hu2018moving}
Y.~Hu, Y.~Fang, Z.~Ge, Z.~Qu, Y.~Zhu, A.~Pradhana, C.~Jiang, A moving least
  squares material point method with displacement discontinuity and two-way
  rigid body coupling, ACM Transactions on Graphics (TOG) 37~(4) (2018) 1--14.

\bibitem{zhang2008material}
D.~Z. Zhang, Q.~Zou, W.~B. VanderHeyden, X.~Ma, Material point method applied
  to multiphase flows, Journal of Computational Physics 227~(6) (2008)
  3159--3173.

\bibitem{homel2017field}
M.~A. Homel, E.~B. Herbold, Field-gradient partitioning for fracture and
  frictional contact in the material point method, International Journal for
  Numerical Methods in Engineering 109~(7) (2017) 1013--1044.

\bibitem{nairn2003material}
J.~A. Nairn, Material point method calculations with explicit cracks, Computer
  Modeling in Engineering and Sciences 4~(6) (2003) 649--664.

\bibitem{homel2018fracture}
M.~Homel, E.~Herbold, Fracture and contact in the material point method: New
  approaches and applications, in: Advances in Computational Coupling and
  Contact Mechanics, World Scientific, 2018, pp. 289--326.

\bibitem{tan2002hierarchical}
H.~Tan, J.~A. Nairn, Hierarchical, adaptive, material point method for dynamic
  energy release rate calculations, Computer Methods in Applied Mechanics and
  Engineering 191~(19-20) (2002) 2123--2137.

\bibitem{zhang2017incompressible}
F.~Zhang, X.~Zhang, K.~Y. Sze, Y.~Lian, Y.~Liu, Incompressible material point
  method for free surface flow, Journal of Computational Physics 330 (2017)
  92--110.

\bibitem{abe2014material}
K.~Abe, K.~Soga, S.~Bandara, Material point method for coupled hydromechanical
  problems, Journal of Geotechnical and Geoenvironmental Engineering 140~(3)
  (2014) 04013033.

\bibitem{guilkey2007eulerian}
J.~Guilkey, T.~Harman, B.~Banerjee, An eulerian--lagrangian approach for
  simulating explosions of energetic devices, Computers \& structures
  85~(11-14) (2007) 660--674.

\bibitem{gaume2018dynamic}
J.~Gaume, T.~Gast, J.~Teran, A.~van Herwijnen, C.~Jiang, Dynamic anticrack
  propagation in snow, Nature communications 9~(1) (2018) 1--10.

\bibitem{yang2017combined}
Y.~Yang, P.~Sun, Z.~Chen, Combined mpm-dem for simulating the interaction
  between solid elements and fluid particles, Communications in Computational
  Physics 21~(5) (2017) 1258--1281.

\bibitem{liu2017multi}
C.~Liu, Q.~Sun, Y.~Yang, Multi-scale modelling of granular pile collapse by
  using material point method and discrete element method, Procedia Engineering
  175 (2017) 29--35.

\bibitem{jiang2020hybrid}
Y.~Jiang, M.~Li, C.~Jiang, F.~Alonso-Marroquin, A hybrid material-point
  spheropolygon-element method for solid and granular material interaction,
  International Journal for Numerical Methods in Engineering 121~(14) (2020)
  3021--3047.

\bibitem{chen2021hybrid}
P.~Y. Chen, M.~Chantharayukhonthorn, Y.~Yue, E.~Grinspun, K.~Kamrin, Hybrid
  discrete-continuum modeling of shear localization in granular media, Journal
  of the Mechanics and Physics of Solids (2021) 104404.

\bibitem{higo2010coupled}
Y.~Higo, F.~Oka, S.~Kimoto, Y.~Morinaka, Y.~Goto, Z.~Chen, A coupled mpm-fdm
  analysis method for multi-phase elasto-plastic soils, Soils and foundations
  50~(4) (2010) 515--532.

\bibitem{raymond2018strategy}
S.~J. Raymond, B.~Jones, J.~R. Williams, A strategy to couple the material
  point method (mpm) and smoothed particle hydrodynamics (sph) computational
  techniques, Computational Particle Mechanics 5~(1) (2018) 49--58.

\bibitem{cummins2002implicit}
S.~Cummins, J.~Brackbill, An implicit particle-in-cell method for granular
  materials, Journal of Computational Physics 180~(2) (2002) 506--548.

\bibitem{guilkey2003implicit}
J.~E. Guilkey, J.~A. Weiss, Implicit time integration for the material point
  method: Quantitative and algorithmic comparisons with the finite element
  method, International Journal for Numerical Methods in Engineering 57~(9)
  (2003) 1323--1338.

\bibitem{homel2016controlling}
M.~A. Homel, R.~M. Brannon, J.~Guilkey, Controlling the onset of numerical
  fracture in parallelized implementations of the material point method (mpm)
  with convective particle domain interpolation (cpdi) domain scaling,
  International Journal for Numerical Methods in Engineering 107~(1) (2016)
  31--48.

\bibitem{goetz2006blast}
D.~Swensen, M.~Denison, J.~Guilkey, T.~Harman, R.~Goetz, A software framework
  for blast event simulation, Report of Reaction Engineering International
  (2006) 9.

\bibitem{lian2011femp}
Y.~Lian, X.~Zhang, X.~Zhou, Z.~Ma, A femp method and its application in
  modeling dynamic response of reinforced concrete subjected to impact loading,
  Computer methods in applied mechanics and engineering 200~(17-20) (2011)
  1659--1670.

\bibitem{jiang2015affine}
C.~Jiang, C.~Schroeder, A.~Selle, J.~Teran, A.~Stomakhin, The affine
  particle-in-cell method, ACM Transactions on Graphics (TOG) 34~(4) (2015)
  1--10.

\bibitem{banerjee2005shell}
B.~Banerjee, Simulation of thin hyperelastic shells with the material point
  method (2005).
\newblock \href {http://arxiv.org/abs/arXiv:cond-mat/0510370}
  {\path{arXiv:arXiv:cond-mat/0510370}}.

\bibitem{zhang2006explicit}
X.~Zhang, K.~Sze, S.~Ma, An explicit material point finite element method for
  hyper-velocity impact, International Journal for Numerical Methods in
  Engineering 66~(4) (2006) 689--706.

\bibitem{lian2012adaptive}
Y.~Lian, X.~Zhang, Y.~Liu, An adaptive finite element material point method and
  its application in extreme deformation problems, Computer methods in applied
  mechanics and engineering 241 (2012) 275--285.

\bibitem{lian2011coupling}
Y.~Lian, X.~Zhang, Y.~Liu, Coupling of finite element method with material
  point method by local multi-mesh contact method, Computer Methods in Applied
  Mechanics and Engineering 200~(47-48) (2011) 3482--3494.

\bibitem{lian2014coupling}
Y.-P. Lian, Y.~Liu, X.~Zhang, Coupling of membrane element with material point
  method for fluid--membrane interaction problems, International Journal of
  Mechanics and Materials in Design 10~(2) (2014) 199--211.

\bibitem{li2021novel}
M.~Li, Y.~Lei, D.~Gao, Y.~Hu, X.~Zhang, A novel material point method (mpm)
  based needle-tissue interaction model, Computer Methods in Biomechanics and
  Biomedical Engineering (2021) 1--15.

\bibitem{cheon2018efficient}
Y.-J. Cheon, H.-G. Kim, An efficient contact algorithm for the interaction of
  material particles with finite elements, Computer Methods in Applied
  Mechanics and Engineering 335 (2018) 631--659.

\bibitem{chen2015improved}
Z.~Chen, X.~Qiu, X.~Zhang, Y.~Lian, Improved coupling of finite element method
  with material point method based on a particle-to-surface contact algorithm,
  Computer Methods in Applied Mechanics and Engineering 293 (2015) 1--19.

\bibitem{wu2018coupled}
B.~Wu, Z.~Chen, X.~Zhang, Y.~Liu, Y.~Lian, Coupled shell-material point method
  for bird strike simulation, Acta Mechanica Solida Sinica 31~(1) (2018) 1--18.

\bibitem{song2020non}
Y.~Song, Y.~Liu, X.~Zhang, A non-penetration fem-mpm contact algorithm for
  complex fluid-structure interaction problems, Computers \& Fluids 213 (2020)
  104749.

\bibitem{bewick2004combined}
B.~T. Bewick, A combined FEM and MPM simulation of impact-resistant design,
  University of Missouri-Columbia, 2004.

\bibitem{aulisa2019monolithic}
E.~Aulisa, G.~Capodaglio, Monolithic coupling of the implicit material point
  method with the finite element method, Computers \& Structures 219 (2019)
  1--15.

\bibitem{larese2019implicit}
A.~Larese, I.~Iaconeta, B.~Chandra, V.~Singer, Implicit mpm and coupled mpm-fem
  in geomechanics, Computational mechanics 175 (2019) 226--232.

\bibitem{chandra2021nonconforming}
B.~Chandra, V.~Singer, T.~Teschemacher, R.~Wuechner, A.~Larese, Nonconforming
  dirichlet boundary conditions in implicit material point method by means of
  penalty augmentation, Acta Geotechnica (2021) 1--21.

\bibitem{guilkey2001implicit}
J.~Guilkey, J.~Weiss, An implicit time integration strategy for use with the
  material point method, in: Proceedings from the First MIT Conference on
  Computational Fluid and Solid Mechanics, Vol.~2, 2001.

\bibitem{nair2012implicit}
A.~Nair, S.~Roy, Implicit time integration in the generalized interpolation
  material point method for finite deformation hyperelasticity, Mechanics of
  Advanced Materials and Structures 19~(6) (2012) 465--473.

\bibitem{charlton2017igimp}
T.~Charlton, W.~Coombs, C.~Augarde, igimp: An implicit generalised
  interpolation material point method for large deformations, Computers \&
  Structures 190 (2017) 108--125.

\bibitem{love2006unconditionally}
E.~Love, D.~L. Sulsky, An unconditionally stable, energy--momentum consistent
  implementation of the material-point method, Computer Methods in Applied
  Mechanics and Engineering 195~(33-36) (2006) 3903--3925.

\bibitem{ortiz1999variational}
M.~Ortiz, L.~Stainier, The variational formulation of viscoplastic constitutive
  updates, Computer methods in applied mechanics and engineering 171~(3-4)
  (1999) 419--444.

\bibitem{kane2000variational}
C.~Kane, J.~E. Marsden, M.~Ortiz, M.~West, Variational integrators and the
  newmark algorithm for conservative and dissipative mechanical systems,
  International Journal for numerical methods in engineering 49~(10) (2000)
  1295--1325.

\bibitem{gast2015optimization}
T.~F. Gast, C.~Schroeder, A.~Stomakhin, C.~Jiang, J.~M. Teran, Optimization
  integrator for large time steps, IEEE Transactions on Visualization and
  Computer Graphics (TVCG) 21~(10) (2015) 1103--1115.

\bibitem{wang2019hierarchical}
X.~Wang, M.~Li, Y.~Fang, X.~Zhang, M.~Gao, M.~Tang, D.~M. Kaufman, C.~Jiang,
  Hierarchical optimization time integration for cfl-rate mpm stepping, ACM
  Transactions on Graphics (TOG) 39~(3) (2020) 1--16.

\bibitem{li2020incremental}
M.~Li, Z.~Ferguson, T.~Schneider, T.~Langlois, D.~Zorin, D.~Panozzo, C.~Jiang,
  D.~M. Kaufman, Incremental potential contact: Intersection-and
  inversion-free, large-deformation dynamics, ACM transactions on graphics
  (2020).

\bibitem{li2020robust}
M.~Li, Robust and accurate simulation of elastodynamics and contact, Ph.D.
  thesis, University of Pennsylvania (2020).

\bibitem{eipc}
Anonymous, Convergent incremental potential contact for piecewise linear
  boundaries, unpublished (N.D.).

\bibitem{Li2021CIPC}
M.~Li, D.~M. Kaufman, C.~Jiang, Codimensional incremental potential contact,
  ACM Trans. Graph. (SIGGRAPH) 40~(4) (2021).

\bibitem{Ferguson:2021:RigidIPC}
Z.~Ferguson, M.~Li, T.~Schneider, F.~Gil-Ureta, T.~Langlois, C.~Jiang,
  D.~Zorin, D.~M. Kaufman, D.~Panozzo, Intersection-free rigid body dynamics,
  ACM Transactions on Graphics (SIGGRAPH) 40~(4) (2021).

\bibitem{Lan2021MIPC}
L.~Lan, Y.~Yang, D.~M. Kaufman, J.~Yao, M.~Li, C.~Jiang, Medial {IPC}:
  Accelerated incremental potential contact with medial elastics, ACM
  Transactions on Graphics (SIGGRAPH) 40~(4) (2021).

\bibitem{choo2021barrier}
J.~Choo, Y.~Zhao, Y.~Jiang, M.~Li, C.~Jiang, K.~Soga, A barrier method for
  frictional contact on embedded interfaces (2021).
\newblock \href {http://arxiv.org/abs/2107.05814} {\path{arXiv:2107.05814}}.

\bibitem{nakamura2021particle}
K.~Nakamura, S.~Matsumura, T.~Mizutani, Particle-to-surface frictional contact
  algorithm for material point method using weighted least squares, Computers
  and Geotechnics 134 (2021) 104069.

\bibitem{tjung2020modeling}
E.~Y. Tjung, S.~Kularathna, K.~Kumar, K.~Soga, Modeling irregular boundaries
  using isoparametric elements in material point method, in: Geo-Congress 2020:
  Modeling, Geomaterials, and Site Characterization, American Society of Civil
  Engineers Reston, VA, 2020, pp. 39--48.

\bibitem{bonet2008nonlinear}
J.~Bonet, R.~D. Wood, Nonlinear Continuum Mechanics for Finite Element
  Analysis, Cambridge University Press, 2008.

\bibitem{radovitzky1999error}
R.~Radovitzky, M.~Ortiz, Error estimation and adaptive meshing in strongly
  nonlinear dynamic problems, Computer Methods in Applied Mechanics and
  Engineering 172~(1-4) (1999) 203--240.

\bibitem{li2019decomposed}
M.~Li, M.~Gao, T.~Langlois, C.~Jiang, D.~M. Kaufman, Decomposed optimization
  time integrator for large-step elastodynamics, ACM Transactions on Graphics
  (TOG) 38~(4) (2019) 1--10.

\bibitem{steffen2008analysis}
M.~Steffen, R.~Kirby, M.~Berzins, Analysis and reduction of quadrature errors
  in the material point method ({MPM}), Int J Numer Meth Eng 76~(6) (2008)
  922--948.

\bibitem{long2019using}
C.~C. Long, G.~Moutsanidis, Y.~Bazilevs, D.~Z. Zhang, Using nurbs as shape
  functions for mpm, Tech. rep., Los Alamos National Lab.(LANL), Los Alamos, NM
  (United States) (2019).

\bibitem{sadeghirad2013second}
A.~Sadeghirad, R.~M. Brannon, J.~Guilkey, Second-order convected particle
  domain interpolation (cpdi2) with enrichment for weak discontinuities at
  material interfaces, International Journal for numerical methods in
  Engineering 95~(11) (2013) 928--952.

\bibitem{nguyen2017family}
V.~P. Nguyen, C.~T. Nguyen, T.~Rabczuk, S.~Natarajan, On a family of convected
  particle domain interpolations in the material point method, Finite Elements
  in Analysis and Design 126 (2017) 50--64.

\bibitem{long2016representing}
C.~Long, D.~Zhang, C.~Bronkhorst, G.~Gray~III, Representing ductile damage with
  the dual domain material point method, Computer Methods in Applied Mechanics
  and Engineering 300 (2016) 611--627.

\bibitem{zhang2013dual}
D.~Z. Zhang, Dual domain material point method for extreme material
  deformation, Tech. rep., Los Alamos National Lab.(LANL), Los Alamos, NM
  (United States) (2013).

\bibitem{jiang2017angular}
C.~Jiang, C.~Schroeder, J.~Teran, An angular momentum conserving
  affine-particle-in-cell method, Journal of Computational Physics 338 (2017)
  137--164.

\bibitem{hammerquist2017new}
C.~C. Hammerquist, J.~A. Nairn, A new method for material point method particle
  updates that reduces noise and enhances stability, Computer methods in
  applied mechanics and engineering 318 (2017) 724--738.

\bibitem{wachter2006implementation}
A.~W{\"a}chter, L.~T. Biegler, On the implementation of an interior-point
  filter line-search algorithm for large-scale nonlinear programming,
  Mathematical programming 106~(1) (2006) 25--57.

\bibitem{brochu2012efficient}
T.~Brochu, E.~Edwards, R.~Bridson, Efficient geometrically exact continuous
  collision detection, ACM Transactions on Graphics (TOG) 31~(4) (2012) 1--7.

\bibitem{moreau2011unilateral}
J.~J. Moreau, On unilateral constraints, friction and plasticity, in: New
  variational techniques in mathematical physics, Springer, 2011, pp. 171--322.

\bibitem{goyal1991planar}
S.~Goyal, A.~Ruina, J.~Papadopoulos, Planar sliding with dry friction part 2.
  dynamics of motion, Wear 143~(2) (1991) 331--352.

\bibitem{teran2005robust}
J.~Teran, E.~Sifakis, G.~Irving, R.~Fedkiw, Robust quasistatic finite elements
  and flesh simulation, in: Proceedings of the 2005 ACM SIGGRAPH/Eurographics
  symposium on Computer animation (SCA), ACM, 2005, pp. 181--190.

\bibitem{chen2008algorithm}
Y.~Chen, T.~A. Davis, W.~W. Hager, S.~Rajamanickam, Algorithm 887: Cholmod,
  supernodal sparse cholesky factorization and update/downdate, ACM
  Transactions on Mathematical Software (TOMS) 35~(3) (2008) 1--14.

\bibitem{nocedal2006numerical}
J.~Nocedal, S.~Wright, Numerical optimization, Springer Science \& Business
  Media, 2006.

\bibitem{li2021lagrangian}
Y.~Li, X.~Li, M.~Li, Y.~Zhu, B.~Zhu, C.~Jiang, Lagrangian--eulerian
  multidensity topology optimization with the material point method,
  International Journal for Numerical Methods in Engineering 122~(14) (2021)
  3400--3424.

\bibitem{huang2010contact}
P.~Huang, X.~Zhang, S.~Ma, X.~Huang,
  \href{https://doi.org/10.1002/nme.2981}{Contact algorithms for the material
  point method in impact and penetration simulation}, International Journal for
  Numerical Methods in Engineering 85~(4) (2010) 498--517.
\newblock \href {https://doi.org/10.1002/nme.2981}
  {\path{doi:10.1002/nme.2981}}.
\newline\urlprefix\url{https://doi.org/10.1002/nme.2981}

\bibitem{barber2009contact}
J.~R. Barber, Plane contact problems, in: Solid Mechanics and Its Applications,
  Springer Netherlands, 2009, pp. 171--198.

\bibitem{liu2020ilsmpm}
C.~Liu, W.~Sun, \href{https://doi.org/10.1016/j.cma.2020.113168}{{ILS}-{MPM}:
  An implicit level-set-based material point method for frictional particulate
  contact mechanics of deformable particles}, Computer Methods in Applied
  Mechanics and Engineering 369 (2020) 113168.
\newblock \href {https://doi.org/10.1016/j.cma.2020.113168}
  {\path{doi:10.1016/j.cma.2020.113168}}.
\newline\urlprefix\url{https://doi.org/10.1016/j.cma.2020.113168}

\bibitem{zimmerman2018febio}
B.~K. Zimmerman, G.~A. Ateshian, \href{https://doi.org/10.1115/1.4040497}{A
  surface-to-surface finite element algorithm for large deformation frictional
  contact in febio}, Journal of Biomechanical Engineering 140~(8) (Jul. 2018).
\newblock \href {https://doi.org/10.1115/1.4040497}
  {\path{doi:10.1115/1.4040497}}.
\newline\urlprefix\url{https://doi.org/10.1115/1.4040497}

\end{thebibliography}

\end{document}